\documentclass{amsart}
\usepackage[T1]{fontenc}
\usepackage{mathtools,amssymb}
\usepackage[hidelinks]{hyperref}
\numberwithin{equation}{section}
\newtheorem{theorem}{Theorem}[section]
\newtheorem{corollary}[theorem]{Corollary}
\newtheorem{lemma}[theorem]{Lemma}
\newtheorem{proposition}[theorem]{Proposition}

\theoremstyle{remark}

\newcommand{\R}{\mathbb R}
\newcommand{\C}{\mathcal C}
\newcommand{\K}{\mathcal K}
\newcommand{\dd}{\,d}
\newcommand{\tr}{\operatorname{tr}}
\allowdisplaybreaks[2]
\usepackage{xcolor}

\title[Uniqueness of a cylindrical cone]{Uniqueness of self-shrinking ends
asymptotic to a cylindrical cone}
\author{Yile Liu}

\address{Department of Mathematics,
Lehigh University,
Bethlehem, PA 18015, USA}

\email{yil924@lehigh.edu}
\date{}
\hypersetup{pdftitle={Uniqueness of self-shrinking ends asymptotic to a cylindircal cone},pdfsubject={Products of regular hypersurface cones with arbitrary Euclidean factors; an explicit compensated Carleman weight}}
\begin{document}
\begin{abstract}

Let $\C^m\subset\R^{m+n}$ be a regular embedded cone, where $n\ge1$, and let $k\ge1$. We prove that two smooth, connected, properly embedded self-shrinking ends whose blow-downs converge locally smoothly to the same cylindrical cone $\C\times\R^k$ away from the axis $\{0\}\times\R^k$ must coincide. As a consequence, we obtain a rigidity result for complete self-shrinkers asymptotic to cylindrical cones. The proof uses Gaussian decay estimates and a sectorial Carleman estimate to establish backward uniqueness.

\end{abstract}
\maketitle
\section{Introduction}
A hypersurface $\Sigma $ in $\R^{m+k+1}$ is called a self-shrinker if
\begin{equation}\label{eq:shrinker}
H=\frac{1}{2}\langle X,\nu\rangle,
\end{equation}
where $X$ is the position vector, $\nu$ is a unit normal, and $ H=-\langle\Delta_gX, \nu\rangle$ is the scalar mean curvature  for the induced metric $g$. More generally, for submanifolds of arbitrary codimension, we denote the mean curvature vector by $\mathbf H=\Delta_gX$, so that $\mathbf H=-H\nu$ in the hypersurface case. The self-shrinker equation for mean curvature vector takes the form
\[
\mathbf H= - \frac{1}{2}X^{\perp}.
\]
For $s<0$, the family $\Sigma_s=\sqrt{-s}\,\Sigma$ is a self-similarly shrinking solution of mean curvature flow. Throughout the paper, we use the reverse-time variable $t=-s>0$, so that $\Sigma_{-t}=\sqrt t\,\Sigma$. Self-shrinkers arise as singularity models of mean curvature flow under parabolic rescaling; see Huisken \cite{Huisken} and Andrews--Chow--Guenther--Langford
\cite[Chapter~13]{ACGL20}. 

The geometry at infinity of noncompact self-shrinkers was investigated by Ilmanen \cite{IlmanenNotes}. He conjectured that a noncompact embedded self-shrinker in $\R^3$ has finitely many ends, each of which is either asymptotically conical or asymptotically cylindrical. More precisely, each end $E$ should either converge locally smoothly to a cone that is smooth away from the origin under blow-down $\lambda E$ as $\lambda\to 0^+$, or converge to a round self-shrinking cylinder of radius $\sqrt{2}$ under translations $E-\tau v$ as $\tau\to\infty$ for some nonzero vector $v$ along its axis. He also conjectured that a self-shrinker with an asymptotically cylindrical end must itself be a cylinder.

Angenent, Ilmanen, and Chopp provided numerical evidence for the existence of noncompact self-shrinkers with two asymptotically conical ends \cite{AIC}. Later, Kapouleas, Kleene, and M{\o}ller \cite{KKM} and X.~H.~Nguyen \cite{NguyenConstruction} constructed complete embedded self-shrinkers with asymptotically conical ends. More recently, Buzano, H.~T.~Nguyen, and Schulz \cite{BNS} constructed noncompact self-shrinkers of arbitrary genus and established asymptotically conical behavior for sufficiently large genus. 
Ketover \cite{Ketover} and Ilmanen and White
\cite{IlmanenWhite} obtained further constructions
of self-shrinkers with two asymptotically conical ends.
Shao and Zou \cite{ShaoZou} constructed embedded
self-shrinkers in $\R^3$ with arbitrarily many
asymptotically conical ends by stacking planes.
These results naturally raise the question whether
two self-shrinking ends asymptotic to the same cone
must coincide.

Wang proved that two smooth, connected, properly embedded self-shrinking ends asymptotic to the same regular cone must coincide \cite{Wang}. She later established an infinite-order uniqueness theorem for asymptotically cylindrical ends \cite{WangCyl}. Bernstein gave a purely elliptic proof of Wang's conical uniqueness theorem using weighted frequency functions \cite{Bernstein}. Khan extended the regular-cone uniqueness theorem to properly immersed self-shrinkers in arbitrary codimension \cite{Khan}. In this paper, we study the corresponding uniqueness problem when the asymptotic cone is of the form $\C\times\R^k$, where $\C$ is a regular cone and $k\ge1$.

Throughout the paper, by a hypersurface we mean a smooth embedded submanifold of codimension one. We write $B_R$ for the open Euclidean ball centered at the origin with radius $R$, and $S_R=\partial B_R$. A hypersurface cone $\C^m\subset\R^{m+1}$ with vertex at the origin is a set of the form
\[
\C=\{r\omega:r\ge0,\ \omega\in\Gamma\},
\]
where $\Gamma$ is a hypersurface of $\mathbb S^m.$  
We say that $\C$ is a regular hypersurface cone if $\Gamma$ is a smooth, connected, closed, embedded submanifold with codimension one. In this case, $\C\setminus\{0\}$ is a smooth embedded hypersurface of $\R^{m+1}$.

More generally, for $n\ge1$, a regular embedded cone $\C^m\subset\R^{m+n}$ is a cone of the same form, where $\Gamma^{m-1}\subset\mathbb S^{m+n-1}$ is a smooth, connected, closed, embedded submanifold. In particular, when $n=1$, this agrees with the definition of a regular hypersurface cone. For $k\ge1$, write
\[
X=(y,z)\in\R^{m+n}\times\R^k,\qquad
r=|y|,\qquad
\rho=|X|=(r^2+|z|^2)^{1/2}.
\]
We denote the axis and the regular part of the cylindrical cone by
\[
\mathcal A=\{0\}\times\R^k,
\qquad
\K=(\C\setminus\{0\})\times\R^k.
\]
The cylindrical cone $\C\times\R^k$ is smooth on $\K$. If $\C$ is nonflat, then $\C\times\R^k$ is singular along $\mathcal A$.

We say that $\lambda M$ converges locally smoothly to $\C\times\R^k$ on $\K$ as $\lambda\to0^+$ if, for every $R>1$ and $j\in\mathbb N$,
\[
\lambda M\cap\{|X|\le R,\ r\ge R^{-1}\}
\]
converges to
\[
(\C\times\R^k)\cap\{|X|\le R,\ r\ge R^{-1}\}
\]
in the $C^j$ topology.

Our first result establishes uniqueness for hypersurface self-shrinking ends.

\begin{theorem}\label{thm:main}
Let $\C^m\subset\R^{m+1}$ be a regular hypersurface cone, let $k\ge1$, and let $R_0>0$. Suppose that $\Sigma$ and $\widetilde\Sigma$ are smooth, connected, properly embedded, two-sided self-shrinking ends in $\R^{m+k+1}\setminus B_{R_0}$, with $\partial\Sigma,\ \partial\widetilde\Sigma\subset S_{R_0}.$
If
\[
\lambda\Sigma\longrightarrow\C\times\R^k,
\qquad
\lambda\widetilde\Sigma\longrightarrow\C\times\R^k
\]
locally smoothly on $\K$ as $\lambda\to0^+$, then $\Sigma=\widetilde\Sigma$.
\end{theorem}

Theorem~\ref{thm:main} extends to properly embedded self-shrinking ends in arbitrary codimension. The proof follows the same strategy, but working with sections of the normal bundle introduces additional difficulties in the graphical equations and their estimates. To highlight the main ideas, we first present the proof in codimension one and discuss the higher-codimension modifications in Section~\ref{extension} and Appendix~\ref{app:higher-codimension}.

\begin{theorem}\label{thm:2}
Let $\C^m\subset\R^{m+n}$ be a regular embedded cone, let $n\ge1$ and $k\ge1$, and let $R_0>0$. Suppose that $M_1$ and $M_2$ are smooth, connected, properly embedded self-shrinking ends in $\R^{m+k+n}\setminus B_{R_0}$, with $\partial M_1,\ \partial M_2\subset S_{R_0}$. If
\[
\lambda M_1\longrightarrow\C\times\R^k,\qquad \lambda M_2\longrightarrow\C\times\R^k
\]
locally smoothly on $\K$ as $\lambda\to0^+$, then $M_1=M_2$.
\end{theorem}

Theorem~\ref{thm:2} has the following consequence for complete self-shrinkers.

\begin{corollary}\label{thm3}
Let $n\ge1$ and $k\ge1$, and let $\C^m\subset\R^{m+n}$ be a regular embedded cone. Suppose that $M_2^m\subset\R^{m+n}$ is a smooth, connected, complete, properly embedded self-shrinker without boundary such that
\[
\lambda M_2\longrightarrow\C
\]
locally smoothly on $\C\setminus\{0\}$ as $\lambda\to0^+$. If $M_1^{m+k}\subset\R^{m+k+n}$ is another smooth, connected, complete, properly embedded self-shrinker without boundary satisfying
\[
\lambda M_1\longrightarrow\C\times\R^k
\]
locally smoothly on $\K$ as $\lambda\to0^+$, then $M_1=M_2\times\R^k$.
\end{corollary}

Related classification results have also been obtained for minimal hypersurfaces asymptotic to cones. Simon and Solomon studied minimal hypersurfaces asymptotic to certain quadratic cones \cite{SimonSolomon}, and Mazet classified minimal hypersurfaces asymptotic to Simons cones up to similarities \cite{Mazet}. More recently, Edelen and Sz\'ekelyhidi proved a Liouville-type theorem for minimal hypersurfaces associated with cylindrical cones under additional geometric and density assumptions \cite{EdelenSzekelyhidi}. Unlike the self-shrinker equation, the minimal hypersurface equation is invariant under dilations, so the asymptotic cone alone does not determine the scale of a nonconical minimal hypersurface.

When $\C$ is nonflat, the main difficulty is the singular axis of $\C\times\R^k$. In this case, Wang's theorem does not apply directly. We therefore carry out the comparison on fixed sectors of the form $\{r>\varepsilon\rho\}$. The normalized annuli in each such sector are compact subsets of $\K$, but no boundary condition is available on $r=\epsilon \rho.$ 

We briefly describe the proof of Theorem~\ref{thm:main}. Following Wang's adaptation \cite{Wang} of the backward-uniqueness method of Escauriaza, Seregin, and \v{S}ver\'ak \cite{ESS}, we write the rescaled ends $\sqrt t\,\Sigma$ and $\sqrt t\,\widetilde\Sigma$ as normal graphs over $\C\times\R^k$ on nested fixed sectors away from the singular axis. 
The local smooth convergence implies that the graph functions, and their difference $w$, converge locally smoothly to zero as $t\downarrow0$. The graphical equations imply that $w$ satisfies a backward parabolic inequality whose lower-order coefficients decay at infinity. A localized estimate of T.~A.~Nguyen \cite{Nguyen}, applied as in Wang's argument, gives Gaussian decay for $w$ and its gradient; this decay is used to remove the cutoff at spatial infinity. The main step is to establish a Carleman estimate of Li--\v{S}ver\'ak type \cite{LS} for $\partial_t+\Delta_{g(t)}$ on $\mathcal D_R=\{X\in\K:\rho>R,\ r>\rho/4\}$. For sufficiently large $a$ and every $u\in C_c^\infty(\mathcal D_R\times(0,1/2))$, we prove
\begin{align*}
&\int_0^{1/2}\!\int_{\mathcal D_R}
e^{2\Phi}\left[|\nabla_gu|^2+\bigl(1+a(\Lambda+\varphi)\bigr)u^2\right]
\dd\mu_g\dd t\\
&\qquad\le C\int_0^{1/2}\!\int_{\mathcal D_R}
e^{2\Phi}|\partial_tu+\Delta_gu|^2\dd\mu_g\dd t,
\end{align*}
where 
\[
\Phi=a\Lambda(t)\varphi+(t+1/2)^2,\qquad
\Lambda(t)=\frac{1/2-t}{(t+1/2)^{3/4}},\qquad
\varphi=r^{3/2}-\frac18\rho^{3/2}.
\] The spatial function $\varphi$ satisfies the positivity estimates required for the Carleman estimate. Moreover, since $\varphi>0$ on $\mathcal{D_R}$, and  $\varphi =0 $ on $r=\rho/4$, we can introduce a cutoff along the level sets of $\Lambda \varphi$ away from this boundary without giving a boundary condition. Combining the Carleman estimate  and Gaussian decay , we obtain $w=0$ on a nonempty open set. The Theorem~\ref{thm:main} follows from elliptic unique continuation.

The higher-codimension extension Theorem~\ref{thm:2} follows the same strategy, with the scalar graph functions replaced by the sections of the normal bundle of $\K$. Their difference $w$ is also a normal bundle section. Calculating the $w$ in fixed ambient Euclidean components, we obtain a vector-valued backward parabolic inequality of the form
\[
|\partial_t w+\Delta_{g(t)}w|
\le C\left(\rho^{-2}|w|+\rho^{-1}|\nabla_{g(t)}w|\right).
\]
The coefficient $\rho^{-1}$ remains uniformly bounded under the parabolic rescaling used for Gaussian decay, while the corresponding first-order term in the cutoff argument can be absorbed by choosing $R$ sufficiently large. Thus the same Gaussian decay and Carleman arguments show that the two ends agree on a nonempty open set. Elliptic unique continuation for the self-shrinker system then proves Theorem~\ref{thm:2}. Applying the result and continuation argument to $M_1$ and $M_2\times\R^k$ gives Corollary~\ref{thm3}.

We expect that the argument can be adapted to properly immersed self-shrinking ends under suitable smooth finite-multiplicity convergence assumptions away from the axis, using the covering-space methods of Khan \cite[Definition~2.3 and Section~3]{Khan}. For simplicity, we restrict the statements and proofs in this paper to the properly embedded case.



The paper is organized as follows. We first prove Theorem~\ref{thm:main} in codimension one. In Section~\ref{sec:geometry}, we obtain the graphical representation, derive the metric estimates, and establish the equation for $w$. Section~\ref{sec:gaussian} proves the Gaussian decay estimate. Section~\ref{sec:carleman} constructs the spatial weight and proves the Carleman estimate for the time-dependent metric. Section~\ref{sec:proof} applies the cutoff argument and completes the proof of Theorem~\ref{thm:main}. In Section~\ref{extension}, we extend the argument to higher codimension and prove Theorem~\ref{thm:2} and Corollary~\ref{thm3}. Appendix~\ref{app:higher-codimension} contains the necessary normal-bundle calculations, the modifications of the scalar estimates, and the elliptic unique-continuation argument.

\section*{Acknowledgments}
I would like to thank my advisors, Ao Sun and Huai-Dong Cao, for their guidance, many helpful discussions, and valuable comments and suggestions on this paper. 

\section*{AI Disclosure}
The author used an LLM as a drafting aid in exploring possible choices of the Carleman weight function.

\section{Geometry of the graphical ends}\label{sec:geometry}

In this section and the next three sections, we work in codimension one. 

Throughout the paper, we use the backward-time variable $t=-s>0$, where $s<0$ is the forward mean curvature flow time. Thus, the time-reversed mean curvature flows are given by $\sqrt t\,\Sigma$ and $\sqrt t\,\widetilde\Sigma$ for $0<t\le1$. Here $m$ and $k$ denote the dimensions of the conical and Euclidean factors, respectively.

We represent these flows as normal graphs over fixed sectors of $\K$ and derive a parabolic differential inequality for the difference of their graph functions. We begin with the geometry of $\K$ on fixed sectors.

On $\K$ let $g_0$ be the induced metric and let $\nu_0=(\nu_\C,0)$ be a fixed unit normal. We write $A_0$ for the second fundamental form and $A_0^\sharp$ for the corresponding shape operator. Thus
\[
(A_0^2)_{ij}=g_0^{pq}(A_0)_{ip}(A_0)_{qj}.
\]
A subscript $0$ refers to $g_0$; a subscript $g$ refers to the induced metric introduced below.

The model metric and its Laplacian are
\begin{align*}
g_0&=dr^2+r^2g_\Gamma+\sum_{\mu=1}^k (dz_\mu)^2,\\
\Delta_0v&=v_{rr}+\frac{m-1}{r}v_r+\frac1{r^2}\Delta_\Gamma v
+\sum_{\mu=1}^k v_{z_\mu z_\mu}.
\end{align*}
The second identity holds for arbitrary $v=v(r,\omega,z)$.

We first record the elementary scaling estimates on fixed sectors.

\begin{lemma}\label{lem:sector}
For each fixed \(0<\varepsilon<1\), the spherical cross-section $\{X\in\K:\rho=1,\ r>\varepsilon\}$ has compact smooth closure with boundary contained in \(\K\). For every
\(L>0\), define the sector $\mathcal A_{L,\varepsilon} = \{X\in\K:L\le\rho\le2L,\ r\ge\varepsilon\rho\}.$ Then \(\mathcal A_{L,\varepsilon}\Subset\K\). Moreover, on
\(\mathcal A_{L,\varepsilon}\),
\[
|\nabla_0^jA_0|\le C_{j,\varepsilon}\rho^{-1-j}, 
\qquad \operatorname{Vol}_{g_0}(\mathcal A_{L,\varepsilon}) \le C_\varepsilon L^{m+k},
\]
for $j= 0, 1, 2,\cdots,$
Here \(C_{j,\varepsilon}=C(m,k,\C,j,\varepsilon)\) and \(C_\varepsilon=C(m,k,\C,\varepsilon)\). 
\end{lemma}

\begin{proof}
Fix \(0<\varepsilon<1\). On \(\mathcal A_{1,\varepsilon}\), we have
\(r\ge\varepsilon\), hence \(\mathcal A_{1,\varepsilon}\Subset\K\).
Therefore, for every integer \(j\ge0\),
\[
\sup_{\mathcal A_{1,\varepsilon}}|\nabla_0^jA_0|
\le C_{j,\varepsilon}.
\]
For \(X\in\mathcal A_{L,\varepsilon}\), write \(X=L\bar X\). Then
\(\bar X\in\mathcal A_{1,\varepsilon}\). By homogeneity,
\[
|\nabla_0^jA_0|(X)
=
L^{-1-j}|\nabla_0^jA_0|(\bar X).
\]
Hence
\[
|\nabla_0^jA_0|(X)
\le C_{j,\varepsilon}L^{-1-j}
\le C_{j,\varepsilon}\rho^{-1-j},
\]
where in the last inequality we enlarged the constant to absorb the factor \(2^{1+j}\), using \(L\le\rho\le2L\).

Finally, since \(\mathcal A_{L,\varepsilon}=L\mathcal A_{1,\varepsilon}\),
\[
\operatorname{Vol}_{g_0}(\mathcal A_{L,\varepsilon})
=
L^{m+k}
\operatorname{Vol}_{g_0}(\mathcal A_{1,\varepsilon})
\le C_\varepsilon L^{m+k}.
\]
\end{proof}

Next, the assumed local smooth convergence gives a uniform graphical representation of the rescaled ends on the fixed sector, as in the normalized-annulus argument of \cite[Section~2]{Wang}.

\begin{lemma}\label{lem:graphs}
There is $R>0$ such that the parts of $\sqrt t\,\Sigma$ and $\sqrt t\,\widetilde\Sigma$ over $\{\rho>R,r>\rho/16\}$ 
are normal graphs
\begin{equation*}
X\longmapsto X+U(X,t)\nu_0(X),\qquad
X\longmapsto X+\widetilde U(X,t)\nu_0(X),\qquad 0<t\le1.
\end{equation*}
The graph functions are smooth for $t>0$ and tend to zero with all local spatial derivatives as $t\downarrow0$. For every fixed finite order $j$,
\begin{equation*}
|\nabla_0^jU|+|\nabla_0^j\widetilde U|\le C_j\rho^{1-j},
\end{equation*}
for $j=0,1,2,\cdots.$
Moreover, $\rho^{-1}(|U|+|\widetilde U|)+|\nabla_0U|+|\nabla_0\widetilde U|$ tends to zero as $\rho\to\infty$, uniformly in $0<t\le1$. Here $C_j=C(m,k,\C,\Sigma,\widetilde\Sigma,j)$.
\end{lemma}
\begin{proof}
Let
\[
\mathcal A
 =\{X\in\K:1/2\le\rho\le2,\ r\ge\rho/16\}
 \Subset\{r>\rho/32\}.
\]
For $\lambda\downarrow0$, write on $\mathcal A$
\[
\lambda\Sigma
 =\{\xi+u_\lambda(\xi)\nu_0(\xi):\xi\in\mathcal A\},
\qquad
\|u_\lambda\|_{C^j(\mathcal A)}\longrightarrow0,
\]
and similarly for $\widetilde\Sigma$. Given $\rho_0\gg1$ and $0<t\le1$, let $\lambda=\frac{\sqrt t}{\rho_0}\le\rho_0^{-1}, \rho_0^{-1}(\sqrt t\,\Sigma)=\lambda\Sigma.$ On $\rho_0\mathcal A$ define $U(\rho_0\xi,t)=\rho_0u_{\sqrt t/\rho_0}(\xi).$ Thus, for $X=\rho_0\xi$,
\[
|\nabla_0^jU(X,t)|
 =\rho_0^{1-j}
 |\nabla_0^ju_{\sqrt t/\rho_0}(\xi)|
 \le C_j\rho_0^{1-j}
 \le C_j\rho^{1-j}.
\]
Moreover,
\[
\sup_{0<t\le1}
 \bigl(\rho^{-1}|U|+|\nabla_0U|\bigr)
 \le C\sup_{0<\lambda\le\rho_0^{-1}}
 \|u_\lambda\|_{C^1(\mathcal A)}
 \longrightarrow0
 \qquad(\rho_0\to\infty).
\]
Uniqueness of normal projection yields equality of the local graph functions on overlapping annuli, hence a single $U$ on $\{\rho>R,r>\rho/16\}$. For every $K'\Subset\{\rho>R,r>\rho/16\}$, $\|\nabla_0^jU(\cdot,t)\|_{C^0(K')}\longrightarrow0 $ as $ t\rightarrow 0.$ The same formulas hold for $\widetilde U$. Smoothness for $t>0$ follows from the smooth normal projection, and $|X+U\nu_0|^2=\rho^2+U^2>tR_0^2.$

\end{proof}

We next compute the induced metric and the equation satisfied by each graph function.
\begin{lemma}
For either graph function $U$, its induced metric is
\begin{equation}\label{eq:metric_formula}
g_U=g_0-2UA_0+U^2A_0^2+dU\otimes dU.
\end{equation}
At each fixed $X\in\K$, the graph function $U$ satisfies the backward mean curvature flow equation

\begin{equation}\label{eq:graph_equation}
\begin{split}
\partial_t U
+ g_U^{ij}\Big[
&(\nabla_0^2 U)_{ij}
    +(A_0)_{ij}
    -U(A_0^2)_{ij} \\
&+\Big((\nabla_i^0 U)(A_0)_{j\ell}
        +(\nabla_j^0 U)(A_0)_{i\ell}
        +U\nabla_i^0(A_0)_{j\ell}
      \Big)
      \Big((I-UA_0^\sharp)^{-1}\nabla_0 U\Big)^\ell
\Big]=0.
\end{split}
\end{equation}
Spatial indices $i,j,\ell,q$ range from $1$ to $m+k$ and are summed.
\end{lemma}
\begin{proof}
Fix $P\in\K$ and choose a $g_0$-geodesic orthonormal frame $\{e_i\}_{i=1}^{m+k}$ at $P$. Let
\[
S_U:=I-UA_0^\sharp,\qquad
W_U:=S_U^{-1}\nabla_0U,
\]

For the normal graph $F=X+U\nu_0,$ write $F_i:=dF(e_i).$ Since $D_{e_i}\nu_0=-(A_0^\sharp)_i{}^j e_j,$ we obtain $F_i=(S_U)_i{}^j e_j+(\nabla_i^0U)\nu_0.$ Therefore
\[
\langle F_i,F_j\rangle
=(g_0)_{ij}-2U(A_0)_{ij}+U^2(A_0^2)_{ij}
+(\nabla_i^0U)(\nabla_j^0U),
\]
which proves \eqref{eq:metric_formula}.

Since $|UA_0^\sharp|=o(1)$, the endomorphism $S_U$ is invertible for sufficiently large $R$. Moreover, $\langle \nu_0-W_U,F_i\rangle=0,$ and hence the unit normal to the graph is $\nu_U=\frac{(\nu_0-W_U)}{\sqrt{1+ |W_U|^2}}.$ We now compute the second covariant derivative of the immersion $F$ with respect to the background connection $\nabla^0$. Define
\[
(\nabla_0^2F)_{ij}
:=
D_{e_i}\bigl(dF(e_j)\bigr)
-
dF\bigl(\nabla_{e_i}^0e_j\bigr).
\]
Because the frame is $g_0$-geodesic at $P$, $\nabla_{e_i}^0e_j=0 $ so at $P$, $(\nabla_0^2F)_{ij}=D_{e_i}F_j.$ At $P$ we also have $D_{e_i}e_j=(A_0)_{ij}\nu_0$ and
\[
D_{e_i}D_{e_j}\nu_0
=
-(\nabla_i^0A_0^\sharp)_j{}^\ell e_\ell
-(A_0^2)_{ij}\nu_0.
\]
Differentiating $F=X+U\nu_0$ twice therefore gives
\begin{align*}
(\nabla_0^2F)_{ij}
={}&
\bigl(
(\nabla_0^2U)_{ij}
+(A_0)_{ij}
-U(A_0^2)_{ij}
\bigr)\nu_0
\\
&-
\bigl(
(\nabla_i^0U)(A_0)_{j\ell}
+(\nabla_j^0U)(A_0)_{i\ell}
+U\nabla_i^0(A_0)_{j\ell}
\bigr)e_\ell .
\end{align*}

Taking the inner product with $\nu_U=\frac{(\nu_0-W_U)}{\sqrt{1+ |W_U|^2}}$ yields
\begin{align*}
\sqrt{1+ |W_U|^2}
\left\langle
(\nabla_0^2F)_{ij},\nu_U
\right\rangle
={}&
(\nabla_0^2U)_{ij}
+(A_0)_{ij}
-U(A_0^2)_{ij}
\\
&+
\bigl(
(\nabla_i^0U)(A_0)_{j\ell}
+(\nabla_j^0U)(A_0)_{i\ell}
+U\nabla_i^0(A_0)_{j\ell}
\bigr)(W_U)^\ell .
\end{align*}

Finally, $F_t=(\partial_tU)\nu_0,$ and hence $\langle F_t,\nu_U\rangle = \frac{\partial_tU}{\sqrt{1+ |W_U|^2}}.$ With our sign for mean curvature, $H_U = -g_U^{ij}\left\langle(\nabla_0^2F)_{ij},\nu_U\right\rangle.$ Since $\sqrt{t}\,\Sigma$ evolves by backward mean curvature flow, $\langle F_t,\nu_U\rangle=H_U.$ Multiplying by $\sqrt{1+ |W_U|^2}$ and substituting the preceding expression gives \eqref{eq:graph_equation}.
\end{proof}

Using the graphical equation and Lemma~\ref{lem:graphs}, we obtain the following sharper estimates.

\begin{lemma}\label{lem:fine}
There is $R>0$ such that, for every fixed finite order $j$, on $\{\rho>R,r>\rho/8\}\times(0,1]$,
\begin{align*}
|\nabla_0^j\partial_tU|+|\nabla_0^j\partial_t\widetilde U|
&\le C_j\rho^{-1-j},\\
|\nabla_0^jU|+|\nabla_0^j\widetilde U|
&\le C_jt\rho^{-1-j}.
\end{align*}
Moreover,
\begin{equation*}
U(X,t)=\sqrt t\,U(X/\sqrt t,1),\qquad
\widetilde U(X,t)=\sqrt t\,\widetilde U(X/\sqrt t,1).
\end{equation*}
Here $C_j=C(m,k,\C,\Sigma,\widetilde\Sigma,j)$.
\end{lemma}
\begin{proof}
For every fixed $q$,
\[
|\nabla_0^qU|\le C_q\rho^{1-q},
\qquad
|\nabla_0^qA_0|\le C_q\rho^{-1-q}.
\]
Let
\[
S_U:=I-UA_0^\sharp,
\qquad
W_U:=S_U^{-1}\nabla_0U,
\]
and define the tensor
\[
(\mathcal T_U)_{ij\ell}
:=(\nabla_i^0U)(A_0)_{j\ell}
 +(\nabla_j^0U)(A_0)_{i\ell}
 +U\nabla_i^0(A_0)_{j\ell}.
\]
Then
\[
|\nabla_0^q(UA_0^\sharp)|
\le C_q\sum_{a+b=q}\rho^{1-a}\rho^{-1-b}
\le C_q\rho^{-q}.
\]
Since $|UA_0^\sharp|=o(1)$,
\[
\nabla_0S_U^{-1}
=S_U^{-1}\nabla_0(UA_0^\sharp)S_U^{-1},
\]
and induction gives
\[
|\nabla_0^qS_U^{-1}|
 +|\nabla_0^qg_U^{-1}|
 +|\nabla_0^qW_U|
\le C_q\rho^{-q}.
\]
Moreover,
\begin{align*}
|\nabla_0^q(\nabla_0^2U+A_0-UA_0^2)|
&\le C_q\rho^{-1-q},\\
|\nabla_0^q\mathcal T_U|
&\le C_q\sum_{a+b=q}
 \bigl(\rho^{-a}\rho^{-1-b}
       +\rho^{1-a}\rho^{-2-b}\bigr)
 \le C_q\rho^{-1-q},\\
|\nabla_0^q(\mathcal T_U\cdot W_U)|
&\le C_q\rho^{-1-q}.
\end{align*}
Applying $\nabla_0^j$ to \eqref{eq:graph_equation} gives
\[
|\nabla_0^j\partial_tU|
\le C_j\rho^{-1-j}.
\]
The same estimate holds for $\widetilde U$. Since
\[
\lim_{s\downarrow0}\nabla_0^jU(X,s)=0,
\qquad
\lim_{s\downarrow0}\nabla_0^j\widetilde U(X,s)=0,
\]
we have
\[
\nabla_0^jU(X,t)
 =\int_0^t\nabla_0^j\partial_sU(X,s)\dd s,
\qquad
|\nabla_0^jU(X,t)|\le C_jt\rho^{-1-j},
\]
and similarly for $\widetilde U$.

For $Y=X/\sqrt t$,
\[
\sqrt t\,[Y+U(Y,1)\nu_0(Y)]
=X+\sqrt t\,U(X/\sqrt t,1)\nu_0(X),
\]
because $\nu_0(\lambda X)=\nu_0(X)$. Uniqueness of normal projection gives
\[
U(X,t)=\sqrt t\,U(X/\sqrt t,1),
\qquad
\widetilde U(X,t)=\sqrt t\,\widetilde U(X/\sqrt t,1).
\]
\end{proof}

Let $w=\widetilde U-U,$ $ g(t)=g_{U(t)}.$ Using $g(t)$ as the reference metric, subtraction of the two graphical equations gives the following estimates.
\begin{proposition}\label{prop:difference}
There are $R,C>0$ such that on $\{\rho>R,r>\rho/8\}\times(0,1]$,
\begin{equation}\label{eq:diff_pde}
\begin{aligned}
\sum_{j=0}^3\rho^j|\nabla_0^j(g(t)-g_0)|&\le Ct\rho^{-2},\\
|\partial_tg(t)|+\rho|\nabla_0\partial_tg(t)|&\le C\rho^{-2},
\qquad \tfrac12g_0\le g(t)\le2g_0,\\
|\partial_tw+\Delta_{g(t)}w|
&\le C\big(\rho^{-2}|w|+t\rho^{-3}|\nabla_{g(t)}w|\big).
\end{aligned}
\end{equation}
The metric norms in the first two lines use $g_0$. For every fixed finite $j$ we also have
\begin{equation*}
|\nabla_0^jw|\le C_jt\rho^{-1-j},\qquad
|\nabla_0^j\partial_tw|\le C_j\rho^{-1-j},
\end{equation*}
and
\begin{equation}\label{eq:w_scaling}
w(X,t)=\sqrt t\,w(X/\sqrt t,1),\qquad
|\nabla_0^jw(X,t)|=t^{(1-j)/2}|\nabla_0^jw(X/\sqrt t,1)|.
\end{equation}
Here $C=C(m,k,\C,\Sigma,\widetilde\Sigma)$ and $C_j=C(m,k,\C,\Sigma,\widetilde\Sigma,j)$.
\end{proposition}
\begin{proof}
The estimates for $w$ and the scaling identity follow directly from Lemma~\ref{lem:fine}. From \eqref{eq:metric_formula},
\[
g-g_0=-2UA_0+U^2A_0^2+dU\otimes dU.
\]
Hence, for $0\le j\le3$,
\[
|\nabla_0^j(g-g_0)|
\le
C_jt\rho^{-2-j}
+C_jt^2\rho^{-4-j}
\le
C_jt\rho^{-2-j}.
\]
Also,
\[
\begin{split}
\partial_tg={}&-2(\partial_tU)A_0
+2U(\partial_tU)A_0^2 
+d(\partial_tU)\otimes dU
+dU\otimes d(\partial_tU),
\end{split}
\]
so
\[
|\partial_tg|+\rho|\nabla_0\partial_tg|
\le C\rho^{-2}.
\]
After increasing $R$,
\[
\frac12g_0\le g\le2g_0.
\]

Write
\[
\widetilde g=g_{\widetilde U},
\qquad
S=I-UA_0^\sharp,
\qquad
\widetilde S=I-\widetilde UA_0^\sharp,
\]
and let
\[
W=S^{-1}\nabla_0U,
\qquad
\widetilde W=\widetilde S^{-1}\nabla_0\widetilde U,
\]
\[
T_{ij\ell}
=
(\nabla_i^0U)(A_0)_{j\ell}
+(\nabla_j^0U)(A_0)_{i\ell}
+U\nabla_i^0(A_0)_{j\ell},
\]
with $\widetilde T$ defined similarly. Then
\[
\widetilde g-g
=
-2wA_0+(\widetilde U+U)wA_0^2
+d\widetilde U\otimes dw+dw\otimes dU,
\]
and therefore
\[
|\widetilde g-g|
+
|\widetilde g^{-1}-g^{-1}|
\le
C\bigl(\rho^{-1}|w|
+t\rho^{-2}|\nabla_0w|\bigr).
\]
Moreover,
\[
\widetilde W-W
=
\widetilde S^{-1}\nabla_0w
+\widetilde S^{-1}wA_0^\sharp S^{-1}\nabla_0U,
\]
hence
\[
|W|+|\widetilde W|
\le Ct\rho^{-2},
\qquad
|\widetilde W-W|
\le
C\bigl(|\nabla_0w|+t\rho^{-3}|w|\bigr).
\]
Similarly,
\[
|T|+|\widetilde T|
\le Ct\rho^{-3},
\qquad
|\widetilde T-T|
\le
C\bigl(\rho^{-1}|\nabla_0w|
+\rho^{-2}|w|\bigr).
\]
The nonlinear terms satisfy
\[
\begin{split}
&\widetilde g^{ij}\widetilde T_{ij\ell}\widetilde W^\ell
-g^{ij}T_{ij\ell}W^\ell\\
&\quad=
(\widetilde g^{ij}-g^{ij})
\widetilde T_{ij\ell}\widetilde W^\ell
+g^{ij}(\widetilde T-T)_{ij\ell}\widetilde W^\ell
+g^{ij}T_{ij\ell}(\widetilde W-W)^\ell,
\end{split}
\]
and consequently
\[
\left|
\widetilde g^{ij}\widetilde T_{ij\ell}\widetilde W^\ell
-g^{ij}T_{ij\ell}W^\ell
\right|
\le
C\bigl(
t\rho^{-4}|w|
+t\rho^{-3}|\nabla_0w|
\bigr).
\]

Subtracting the two equations \eqref{eq:graph_equation} gives
\[
\partial_tw+g^{ij}(\nabla_0^2w)_{ij}
=
E_1+E_2+E_3+E_4,
\]
where
\[
E_1
=
-(\widetilde g^{ij}-g^{ij})
(\nabla_0^2\widetilde U)_{ij},
\]
\[
E_2
=
-(\widetilde g^{ij}-g^{ij})(A_0)_{ij},
\]
\[
E_3
=
(\widetilde g^{ij}\widetilde U-g^{ij}U)(A_0^2)_{ij},
\]
and
\[
E_4
=
-\left(
\widetilde g^{ij}\widetilde T_{ij\ell}\widetilde W^\ell
-g^{ij}T_{ij\ell}W^\ell
\right).
\]
Using
\[
|\nabla_0^2\widetilde U|\le Ct\rho^{-3},
\qquad
|A_0|\le C\rho^{-1},
\qquad
|A_0^2|\le C\rho^{-2},
\]
we obtain
\[
|E_1|
\le
C\bigl(
t\rho^{-4}|w|
+t^2\rho^{-5}|\nabla_0w|
\bigr),
\]
\[
|E_2|
\le
C\bigl(
\rho^{-2}|w|
+t\rho^{-3}|\nabla_0w|
\bigr),
\]
and, since
\[
\widetilde g^{ij}\widetilde U-g^{ij}U
=
g^{ij}w+(\widetilde g^{ij}-g^{ij})\widetilde U,
\]
\[
|E_3|
\le
C\bigl(
\rho^{-2}|w|
+t^2\rho^{-5}|\nabla_0w|
\bigr).
\]
Together with the estimate for $E_4$,
\[
\left|
\partial_tw+g^{ij}(\nabla_0^2w)_{ij}
\right|
\le
C\bigl(
\rho^{-2}|w|
+t\rho^{-3}|\nabla_0w|
\bigr).
\]

Finally,
\[
(\Gamma_g-\Gamma_0)^\ell_{ij}
=
\frac12g^{\ell q}
\bigl(
\nabla_i^0g_{jq}
+\nabla_j^0g_{iq}
-\nabla_q^0g_{ij}
\bigr),
\]
so
\[
|\Gamma_g-\Gamma_0|
\le Ct\rho^{-3}.
\]
Since
\[
\Delta_gw
=
g^{ij}(\nabla_0^2w)_{ij}
-g^{ij}(\Gamma_g-\Gamma_0)^\ell_{ij}\nabla_\ell^0w,
\]
we conclude that
\[
|\partial_tw+\Delta_gw|
\le
C\bigl(
\rho^{-2}|w|
+t\rho^{-3}|\nabla_0w|
\bigr).
\]
Using $\frac12g_0\le g\le2g_0$ gives
\[
|\partial_tw+\Delta_gw|
\le
C\bigl(
\rho^{-2}|w|
+t\rho^{-3}|\nabla_gw|
\bigr),
\]
which proves \eqref{eq:diff_pde}.
\end{proof}

\section{Gaussian decay}\label{sec:gaussian}
Throughout Sections~\ref{sec:gaussian} and~\ref{sec:carleman}, $C\geq1$ and $c\in(0,1)$ denote constants that may change from line to line but depending on $m,k, \mathcal{C}, \Sigma, \tilde{\Sigma}$. 

By Proposition \ref{prop:difference}, the function $w$ and its spatial gradient vanish at $t=0.$ In this section, we establish Gaussian decay estimates for $w$ and its gradient away from the singular axis.

\begin{proposition}\label{lem:gaussian}
There are $R,C,c>0$ such that
\begin{equation}\label{eq:gaussian}
|w(X,t)|+|\nabla_{g(t)}w(X,t)|
\le C\exp(-c\rho^2/t)
\end{equation}
for $\rho>R$, $r\ge\rho/4$, and $0<t\le1$.
Here $C=C(m,k,\C,\Sigma,\widetilde\Sigma)$ and $c=c(m,k,\C,\Sigma,\widetilde\Sigma)$.
\end{proposition}


To prove Proposition~\ref{lem:gaussian}, we first establish a local estimate for backward parabolic equations with vanishing initial condition. We use the Carleman estimate of Nguyen\cite[Lemma 2.1]{Nguyen}, together with the localization argument in Wang\cite[Lemma 3.6]{Wang}.

Let $d\ge1$. Let $a^{ij}$ be symmetric on $B_4\times[0,T]\subset\R^d\times[0,T]$ and suppose
\[
\lambda_{1} I\le(a^{ij})\le\lambda_{1}^{-1}I,
\qquad a(0,0)=I,
\qquad
|a(x,t)-a(y,s)|\le L_0\bigl(|x-y|+|t-s|^{1/2}\bigr),
\]
where $0<\lambda_{1}\le1$. Assume also that $a^{ij}$ is smooth for positive time and has bounded first spatial derivatives. There is $N\ge1$, depending on $d$, $\lambda_1$, and $L_0$, such that for every $\alpha\ge2$ there is a positive function $\sigma$ on $(0,4/\alpha)$ satisfying $N^{-1}s\le\sigma(s)\le s.$ For $0<\beta<1/\alpha$, every compactly supported $v$ on $\R^d\times[0,2/\alpha)$ with $v(\cdot,0)=0$ satisfies
\begin{equation}\label{eq:local_carleman}
\begin{split}
&\int\bigl(\alpha^2v^2+\alpha\sigma(t+\beta)|Dv|^2\bigr)
\sigma(t+\beta)^{-\alpha}G_\beta\dd x\dd t\\
&\quad\le C\int|\partial_tv+\partial_i(a^{ij}\partial_jv)|^2
\sigma(t+\beta)^{1-\alpha}G_\beta\dd x\dd t
+C^\alpha\alpha^\alpha\sup_t\int(v^2+|Dv|^2)\dd x,
\end{split}
\end{equation}
where
\[
G_\beta(x,t)=(t+\beta)^{-d/2}
\exp\left(-\frac{|x|^2}{4(t+\beta)}\right).
\]

\begin{lemma}\label{lem:local_l2}
Suppose
\[
|\partial_tu+\partial_i(a^{ij}\partial_ju)|
\le K(|u|+|Du|),
\qquad
|u|+|Du|\le M,
\]
and
\[
u(\cdot,0)=Du(\cdot,0)=0.
\]
There are $\delta_*\in(0,1/4)$, $c_*\in(0,1)$, and $c,C,\alpha_0>0$ such that, whenever $\alpha\ge\alpha_0$ and $4/\alpha\le T$,
\begin{equation}\label{eq:local_l2_claim}
\int_{B_{2\delta_*}\times[0,c_*\delta_*^2/\alpha]}
(u^2+|Du|^2)\dd x\dd t
\le CM^2e^{-c\alpha}.
\end{equation}
The same conclusion holds if $u$ is smooth only for positive time, provided $u=O(t)$ locally and the required derivatives remain bounded up to $t=0$. The constants depend on $d$, $\lambda_1$, $L_0$, and $K$.
\end{lemma}

\begin{proof}
Let $\eta\in C_c^\infty(B_2)$ satisfy $\eta=1$ on $B_1$, and let $\chi_\alpha$ be a time cutoff supported in $[0,2/\alpha]$ and equal to one on $[0,1/\alpha]$. Let
\[
\widehat u=\eta\chi_\alpha u.
\]
After extending $a^{ij}$ outside $B_3$ with the same coefficient bounds,
\[
|\partial_t\widehat u+\partial_i(a^{ij}\partial_j\widehat u)|
\le C(|\widehat u|+|D\widehat u|)+C\alpha M\mathbf1_E,
\]
where $E$ is contained in the cutoff regions. For large $\alpha$, the first term is absorbed in \eqref{eq:local_carleman}, while the cutoff terms give
\[
\int\bigl(\alpha^2\widehat u^2+\alpha\sigma|D\widehat u|^2\bigr)
\sigma^{-\alpha}G_\beta
\le CM^2C_0^\alpha\alpha^{\alpha+d/2}.
\]

Choose $\delta_*>0$ small and let
\[
\beta=\frac{c_*\delta_*^2}{2\alpha},
\qquad
Q_{\delta_*}=B_{2\delta_*}\times
\left[0,\frac{c_*\delta_*^2}{\alpha}\right],
\]
where $c_*>0$ is fixed sufficiently small. Then $\widehat u=u$ on $Q_{\delta_*}$, and on this cylinder
\[
c\frac{\delta_*^2}{\alpha}
\le\sigma(t+\beta)\le
C\frac{\delta_*^2}{\alpha},
\qquad
G_\beta(x,t)
\ge C^{-1}\alpha^{d/2}\delta_*^{-d}e^{-C_1\alpha}.
\]
Hence both weights on the left of \eqref{eq:local_carleman} are bounded from below, up to fixed factors, by
\[
\alpha^{\alpha+d/2}
\left(c_1\delta_*^{-2}e^{-C_1}\right)^\alpha.
\]
It follows that
\[
\int_{Q_{\delta_*}}(u^2+|Du|^2)\dd x\dd t
\le
CM^2\left(C_2e^{C_1}\delta_*^2\right)^\alpha.
\]
Choose $\delta_*$ so that $C_2e^{C_1}\delta_*^2<1$. Then
\[
\int_{Q_{\delta_*}}(u^2+|Du|^2)\dd x\dd t
\le CM^2e^{-c\alpha},
\]
which proves \eqref{eq:local_l2_claim}.

For the weaker initial regularity, insert a time cutoff near $t=0$ and let its support shrink to zero. The additional weighted error tends to zero by the assumed bounds, and the same estimate follows.
\end{proof}

We now apply the preceding estimate to the rescaled equation for $w$ in tangent-plane coordinates and use its scaling law to obtain Gaussian decay.

\begin{proof}[Proof of Proposition~\ref{lem:gaussian}]
We follow the local argument of \cite[Lemma~3.6]{Wang}. Fix $X_0\in\K$ with $r(X_0)\ge\rho_0/4$, where $\rho_0=|X_0|$. If $|X-X_0|\le\kappa\rho_0$, then $r(X)/\rho(X)\ge(1/4-\kappa)/(1+\kappa)$. Choose $\kappa>0$ small enough that this is larger than $1/8$.


By compactness of the normalized sector $\{ \rho=1, r\geq 1/8\}\subseteq \K$ and the homogeneity of $\mathcal{K}$, there exists $c_{0}>0,$ independent of $X_0$, such that $\K$ can be represented as a graph over $T_{X_0}\K$ in a neighborhood of radius $4\ell$, for $\ell:=c_{0}\rho_0$. Taking $c_0$ sufficiently small, this neighborhood is contained in $\{r> \rho/8\}.$

In these coordinates, let $a^{ij}=g^{ij}$. Then
\[
C^{-1}I\le(a^{ij})\le CI,\qquad
|Da|\le C\rho_0^{-1},\qquad
|\partial_ta|\le C\rho_0^{-2},
\]
and, after choosing orthonormal coordinates at $(X_0,0)$, $a(0,0)=I$. Moreover,
\[
\Delta_gw
=\partial_i(a^{ij}\partial_jw)
+a^{ij}\partial_i\log\sqrt{\det g}\,\partial_jw,
\qquad
\left|D\log\sqrt{\det g}\right|\le C\rho_0^{-1}.
\]
Consequently, Proposition \ref{prop:difference} gives
\[
|\partial_tw+\partial_i(a^{ij}\partial_jw)|
\le C\bigl(\rho_0^{-2}|w|+\rho_0^{-1}|Dw|\bigr).
\]

Let $\Psi$ denote the tangent-plane chart and let
\[
u(x,\tau)=w(\Psi(\ell x),\ell^2\tau),\qquad
b^{ij}(x,\tau)=a^{ij}(\Psi(\ell x),\ell^2\tau).
\]
Then on $B_4\times[0,\ell^{-2}]$,
\[
|\partial_\tau u+\partial_i(b^{ij}\partial_ju)|
\le C(|u|+|Du|),
\qquad
|u|+|Du|\le C\rho_0^{-1},
\]
and the coefficients $b^{ij}$ satisfy the hypotheses of Lemma~\ref{lem:local_l2} with uniform constants. Moreover, $u(\cdot,0)=Du(\cdot,0)=0$.

Choose $A$ sufficiently large and let $\alpha=A\ell^2$. Lemma~\ref{lem:local_l2} gives
\[
\int_{B_{2\delta_*}\times[0,c_*\delta_*^2/\alpha]}
(u^2+|Du|^2)\,dx\,d\tau
\le C\rho_0^{-2}e^{-c_1\alpha}.
\]
The derivative estimates in Proposition \ref{prop:difference} give
\[
|D^2u|\le C\rho_0^{-1},\qquad
|\partial_\tau u|+|\partial_\tau Du|\le C\rho_0.
\]

Choose $0<\gamma<c_1/(2(m+k+1))$ and let $\tau_*=\frac{c_*\delta_*^2}{2\alpha},$ $r_\alpha=\frac{c_*\delta_*^2}{8\alpha}e^{-\gamma\alpha}.$
For large $\alpha$,
\[
B_{r_\alpha}^{m+k+1}(0,\tau_*)
\subset
B_{2\delta_*}\times(0,c_*\delta_*^2/\alpha).
\]
For each component $f$ of $(u,Du)$,
\[
\begin{aligned}
|f(0,\tau_*)|
&\le
Cr_\alpha^{-(m+k+1)/2}
\|f\|_{L^2}
+C\rho_0r_\alpha\\
&\le
C\rho_0^{-1}
\alpha^{(m+k+1)/2}
e^{-c\alpha}.
\end{aligned}
\]
Since $\alpha=A\ell^2$, after decreasing $c$,
\[
|w(X_0,t_*)|+|\nabla_0w(X_0,t_*)|
\le Ce^{-c\rho_0^2},
\qquad
t_*=\ell^2\tau_*=\frac{c_*\delta_*^2}{2A}.
\]
The time $t_*>0$ is independent of $X_0$. Using \eqref{eq:w_scaling} with $X_0=\sqrt{t_*}\,X$ gives
\[
|w(X,1)|+|\nabla_0w(X,1)|
\le Ce^{-c\rho^2}.
\]
Applying \eqref{eq:w_scaling} once more,
\[
|w(X,t)|
\le C\sqrt t\,e^{-c\rho^2/t},
\qquad
|\nabla_0w(X,t)|
\le Ce^{-c\rho^2/t}.
\]
Since $\frac12g_0\le g(t)\le2g_0$, \eqref{eq:gaussian} follows.
\end{proof}

\section{A Carleman estimate}\label{sec:carleman}
In this section we prove the Carleman estimate used in the backward-uniqueness argument. We construct an exponential Carleman weight involving a spatial function $\varphi$ and a time-dependent factor. Since no boundary condition is available on $r=\rho/4$, we choose $\varphi$ to vanish there, so that the cutoff can be supported away from this boundary while retaining the positivity estimates needed for the Carleman estimate.
We fix
\begin{equation}\label{eq:phase_pair}
\varphi=r^{3/2}-\frac18\rho^{3/2},
\qquad
\widetilde Q=\frac3{16}\rho^{-1/2},
\qquad
Q=\widetilde Q+\delta_0\rho^{-1/2},
\end{equation}
where $\delta_0>0$ will be chosen below, and let
\[
\mathcal D_R=\{X\in\K:\rho>R,\ r>\rho/4\}.
\]
Then $\varphi>0$ on $\mathcal D_R$ and $\varphi=0$ on $r=\rho/4$.

\begin{lemma}\label{lem:hessian}
The spatial function $\varphi$ in \eqref{eq:phase_pair} satisfies
\[
\nabla_0\varphi
=\frac32\left(r^{1/2}-\frac18\rho^{-1/2}r\right)\partial_r
-\frac3{16}\rho^{-1/2}\sum_{\mu=1}^kz_\mu\partial_{z_\mu},
\]
and
\[
\begin{aligned}
\nabla_0^2\varphi
={}&\frac34r^{-1/2}dr\otimes dr
+\frac32r^{-1/2}(r^2g_\Gamma)
-\frac3{16}\rho^{-1/2}g_0
+\frac3{32}\rho^{-1/2}d\rho\otimes d\rho.
\end{aligned}
\]
Consequently,
\begin{equation}\label{eq:compensated_hessian}
\begin{aligned}
\nabla_0^2\varphi+\widetilde Qg_0
={}&\frac34r^{-1/2}dr\otimes dr
+\frac32r^{-1/2}(r^2g_\Gamma)
+\frac3{32}\rho^{-1/2}d\rho\otimes d\rho
\ge0.
\end{aligned}
\end{equation}
\end{lemma}

\begin{proof}
Fix the point under consideration and choose $\{\bar e_a\}$ to be a $g_\Gamma$-orthonormal frame geodesic at that point. Let $e_a=r^{-1}\bar e_a$. Then, at this point,
\[
\nabla^0_{\partial_r}\partial_r=\nabla^0_{\partial_r}e_a=0,
\qquad
\nabla^0_{e_a}\partial_r=r^{-1}e_a,
\qquad
\nabla^0_{e_a}e_b=-r^{-1}\delta_{ab}\partial_r.
\]
All connection terms involving $\partial_{z_\mu}$ vanish. Hence
\[
\nabla_0^2r^{3/2}
=\frac34r^{-1/2}dr\otimes dr
+\frac32r^{-1/2}(r^2g_\Gamma),
\]
while
\[
\nabla_0^2\rho^{3/2}
=\frac32\rho^{-1/2}g_0
-\frac34\rho^{-1/2}d\rho\otimes d\rho.
\]
The formulas follow by differentiation and subtraction.
\end{proof}

The negative scalar term in the Hessian is canceled by $\widetilde Qg_0$, as in the scalar compensation of \cite[Proposition~2.3 and Section~4]{LS}. We add a small positive term to obtain strict tensor positivity while preserving the remaining positive margins.

\begin{lemma}\label{lem:strict_weight}
For a sufficiently small $\delta_0>0$, there are $C,c>0$ such that on $r\ge\rho/4$,
\begin{equation}\label{eq:strict_weight}
\begin{aligned}
c\rho\le|\nabla_0\varphi|^2&\le C\rho,\\
\nabla_0^2\varphi+Qg_0&\ge c\rho^{-1/2}g_0,\\
\nabla_0^2\varphi(\nabla_0\varphi,\nabla_0\varphi)
-Q|\nabla_0\varphi|^2&\ge c\rho^{1/2},\\
\frac12|\nabla_0\varphi|^2-\varphi Q&\ge c\rho,\\
|\Delta_0^2\varphi|+|\Delta_0Q|&\le C\rho^{-5/2}.
\end{aligned}
\end{equation}
The hypersurface $r=\rho/4$ is a regular level set of $\varphi$. Here $\delta_0=\delta_0(m,k)$, $C=C(m,k)$, and $c=c(m,k)$.
\end{lemma}
\begin{proof}
Let $s=r/\rho\in[1/4,1]$. Since
\[
\partial_r\varphi
=\frac32\rho^{1/2}\left(s^{1/2}-\frac18s\right),
\]
homogeneity gives
\[
c\rho\le|\nabla_0\varphi|^2\le C\rho.
\]
The second inequality follows from \eqref{eq:compensated_hessian} and $Q-\widetilde Q=\delta_0\rho^{-1/2}$.

Put $x=s^{1/2}\in[1/2,1]$. Direct substitution gives
\[
\nabla_0^2\varphi(\nabla_0\varphi,\nabla_0\varphi)
-\widetilde Q|\nabla_0\varphi|^2
=c_0\rho^{1/2}h(x),
\]
where
\[
h(x)=64x^6+56x^3-384x^2+512x-3
=64x^6+56x^3+128x(4-3x)-3.
\]
Thus $h(x)\ge c_1>0$ on $[1/2,1]$. Replacing $\widetilde Q$ by $Q$ subtracts at most $C\delta_0\rho^{1/2}$, and the third inequality follows for small $\delta_0$.

Similarly,
\[
\frac12|\nabla_0\varphi|^2-\varphi\widetilde Q
=c_2\rho\bigl(192x^2-80x^3+7\bigr)
\ge c_3\rho.
\]
Since $0\le\varphi\le\rho^{3/2}$, the additional term in $Q$ is bounded by $C\delta_0\rho$, which proves the fourth inequality.

The functions $\varphi$ and $Q$ are homogeneous of degrees $3/2$ and $-1/2$. Hence
\[
|\Delta_0^2\varphi|+|\Delta_0Q|\le C\rho^{-5/2}
\]
on the closed sector. Finally,
\[
\varphi=\rho^{3/2}\left(s^{3/2}-\frac18\right),
\qquad
\left.\frac{d}{ds}s^{3/2}\right|_{s=1/4}\ne0,
\]
so $r=\rho/4$ is a regular level set.
\end{proof}

The preceding estimates are for $g_0$. We next transfer them to the induced metric $g(t)$.

\begin{lemma}\label{lem:transfer}
There are $R,C,c>0$ such that the functions in \eqref{eq:phase_pair} satisfy on $\mathcal D_R\times(0,1]$
\begin{equation}\label{eq:transfer_estimates}
\begin{aligned}
c\rho\le|\nabla_g\varphi|^2&\le C\rho,\\
\nabla_g^2\varphi+Qg&\ge c\rho^{-1/2}g,\\
\nabla_g^2\varphi(\nabla_g\varphi,\nabla_g\varphi)
-Q|\nabla_g\varphi|^2&\ge c\rho^{1/2},\\
\frac12|\nabla_g\varphi|^2-\varphi Q&\ge c\rho,\\
|\Delta_g^2\varphi|+|\Delta_gQ|&\le C\rho^{-5/2}.
\end{aligned}
\end{equation}
Here $C=C(m,k,\C,\Sigma)$ and $c=c(m,k,\C,\Sigma)$.
\end{lemma}
\begin{proof}
By Proposition \ref{prop:difference},
\[
|\nabla_0^j(g^{-1}-g_0^{-1})|\le Ct\rho^{-2-j},
\quad 0\le j\le3,
\]
and
\[
|\nabla_0^j(\Gamma_g-\Gamma_0)|\le Ct\rho^{-3-j},
\quad 0\le j\le2.
\]
Since $|\nabla_0^j\varphi|\le C_j\rho^{3/2-j}$,
\[
|\nabla_g\varphi-\nabla_0\varphi|\le Ct\rho^{-3/2},
\qquad
|\nabla_g^2\varphi-\nabla_0^2\varphi|\le Ct\rho^{-5/2},
\]
\[
\big||\nabla_g\varphi|^2-|\nabla_0\varphi|^2\big|
\le Ct\rho^{-1},
\qquad
|Q(g-g_0)|\le Ct\rho^{-5/2},
\]
and
\[
\begin{aligned}
&\big|
\nabla_g^2\varphi(\nabla_g\varphi,\nabla_g\varphi)
-\nabla_0^2\varphi(\nabla_0\varphi,\nabla_0\varphi)
\big| 
\le Ct\rho^{-3/2}.
\end{aligned}
\]
The first four inequalities in \eqref{eq:transfer_estimates} now follow from Lemma~\ref{lem:strict_weight} after increasing $R$.

For any smooth $v$,
\[
(\Delta_g-\Delta_0)v
=(g^{ij}-g_0^{ij})(\nabla_0^2v)_{ij}
-g^{ij}(\Gamma_g-\Gamma_0)^\ell_{ij}\nabla_\ell^0v.
\]
Thus
\[
|\nabla_0^j((\Delta_g-\Delta_0)\varphi)|
\le Ct\rho^{-5/2-j},
\quad 0\le j\le2.
\]
Since $\Delta_0\varphi$ and $Q$ are homogeneous of degree $-1/2$,
\[
|(\Delta_g-\Delta_0)(\Delta_0\varphi)|
+|(\Delta_g-\Delta_0)Q|
\le Ct\rho^{-9/2}.
\]
Using
\[
\Delta_g^2\varphi-\Delta_0^2\varphi
=(\Delta_g-\Delta_0)(\Delta_0\varphi)
+\Delta_g((\Delta_g-\Delta_0)\varphi),
\]
we obtain
\[
|\Delta_g^2\varphi-\Delta_0^2\varphi|
+|\Delta_gQ-\Delta_0Q|
\le Ct\rho^{-9/2}.
\]
The final inequality follows from Lemma~\ref{lem:strict_weight}.
\end{proof}

We next choose the time factor. For $0<t<1/2$, let
\begin{equation}\label{eq:spacetime_weight}
\Lambda(t)=\frac{1/2-t}{(t+1/2)^{3/4}},
\qquad
\Phi(X,t)=a\Lambda(t)\varphi(X)+(t+1/2)^2,
\qquad
F=4a\Lambda Q+1.
\end{equation}
This is the shifted factor from \cite[Proposition~2.3]{LS}.

\begin{lemma}\label{lem:time}
There are $c,C>0$ such that on $0<t<1/2$,
\begin{equation}\label{eq:time_bounds}
-\Lambda'\ge c,
\qquad
0\le\Lambda\le C(-\Lambda'),
\qquad
1-\frac{\Lambda}{2(-\Lambda')}\ge\frac12+c,
\end{equation}
and
\[
\Lambda''+\Lambda'\ge c,
\qquad
\Lambda''+\Lambda'\ge c(-\Lambda').
\]
The constants $c$ and $C$ are absolute.
\end{lemma}
\begin{proof}
Let $s=t+1/2\in[1/2,1)$. Then
\[
-\Lambda'=s^{-7/4}\left(\frac34+\frac14s\right),
\]
\[
\Lambda''+\Lambda'
=s^{-11/4}\left(\frac{21}{16}-\frac9{16}s-\frac14s^2\right),
\qquad
\frac{\Lambda}{-\Lambda'}
=\frac{s(1-s)}{\frac34+\frac14s}.
\]
The bounds follow on $[1/2,1]$.
\end{proof}

We now combine the spatial and time estimates. Following \cite[Proposition~2.3]{LS}, we apply the usual square-completion argument to the conjugated equation. Since the metric depends on time, we keep the resulting $\partial_tg$ terms explicit.

\begin{proposition}\label{prop:carleman}
There are $R,a_0,C>0$ such that, for every $a\ge a_0$ and every $u\in C_c^\infty(\mathcal D_R\times(0,1/2))$,
\begin{equation}\label{eq:carleman}
\begin{split}
&\int_0^{1/2}\!\int_{\mathcal D_R}
e^{2\Phi}\left[|\nabla_gu|^2+\bigl(1+a(\Lambda+\varphi)\bigr)u^2\right]
\dd\mu_g\dd t\\
&\hspace{12mm}\le C\int_0^{1/2}\!\int_{\mathcal D_R}
e^{2\Phi}|\partial_tu+\Delta_gu|^2\dd\mu_g\dd t.
\end{split}
\end{equation}
Here $C=C(m,k,\C,\Sigma)$.
\end{proposition}
\begin{proof}
All integrals are over $\mathcal D_R\times(0,1/2)$ with respect to $\dd\mu_g\dd t$. Let $Y=e^\Phi u$ and $V=|\nabla_g\Phi|^2-\Phi_t$. Then
\[
e^\Phi(\partial_tu+\Delta_gu)
=
[\Delta_gY+VY]
+
[Y_t-2\langle\nabla_g\Phi,\nabla_gY\rangle-(\Delta_g\Phi)Y].
\]

Moreover, $\partial_t\dd\mu_g=\frac12\tr_g(\partial_tg)\dd\mu_g$ and
\[
\partial_t|\nabla_gY|^2
=
-(\partial_tg)(\nabla_gY,\nabla_gY)
+2\langle\nabla_gY,\nabla_gY_t\rangle.
\]
Integration by parts gives
\[
2\int(\Delta_gY)Y_t
=
\int
\left[
\frac12\tr_g(\partial_tg)g-\partial_tg
\right](\nabla_gY,\nabla_gY),
\]
and
\begin{align*}
-4\int\Delta_gY\langle\nabla_g\Phi,\nabla_gY\rangle
&=
4\int\nabla_g^2\Phi(\nabla_gY,\nabla_gY)
-2\int(\Delta_g\Phi)|\nabla_gY|^2,\\
-2\int\Delta_gY(\Delta_g\Phi)Y
&=
2\int(\Delta_g\Phi)|\nabla_gY|^2
-\int(\Delta_g^2\Phi)Y^2.
\end{align*}
The remaining cross terms are
\[
\int
\left[
-V_t-\frac12\tr_g(\partial_tg)V
+2\langle\nabla_g\Phi,\nabla_gV\rangle
\right]Y^2.
\]

Using $0\le\int(\Delta_gY+VY+FY/2)^2$, we obtain
\[
\int(\Delta_gY+VY)^2
\ge
\int F|\nabla_gY|^2
+
\int
\left[
-\frac12\Delta_gF-VF-\frac14F^2
\right]Y^2.
\]

We first consider the terms containing $\nabla_gY$. Since $\nabla_g^2\Phi=a\Lambda\nabla_g^2\varphi$ and $F=4a\Lambda Q+1$, their coefficient is
\[
\mathcal G
=
4a\Lambda(\nabla_g^2\varphi+Qg)
+
\left(1+\frac12\tr_g(\partial_tg)\right)g
-\partial_tg.
\]
By \eqref{eq:transfer_estimates} and $|\partial_tg|+|\tr_g(\partial_tg)|\le C\rho^{-2}$, after increasing $R$ we have $\mathcal G\ge cg$.

It remains to identify the zero-order terms. Let $s=t+1/2$. Since $V=a^2\Lambda^2|\nabla_g\varphi|^2-a\Lambda'\varphi-2s$, we have
\begin{align*}
-V_t
={}&
-2a^2\Lambda\Lambda'|\nabla_g\varphi|^2
+a^2\Lambda^2(\partial_tg)(\nabla_g\varphi,\nabla_g\varphi)
+a\Lambda''\varphi+2,\\
2\langle\nabla_g\Phi,\nabla_gV\rangle
={}&
4a^3\Lambda^3
\nabla_g^2\varphi(\nabla_g\varphi,\nabla_g\varphi)
-2a^2\Lambda\Lambda'|\nabla_g\varphi|^2,\\
-FV
={}&
-4a^3\Lambda^3Q|\nabla_g\varphi|^2
-a^2\Lambda^2|\nabla_g\varphi|^2
+4a^2\Lambda\Lambda'\varphi Q\\
&\quad
+a\Lambda'\varphi+8as\Lambda Q+2s,\\
-\frac14F^2
={}&
-4a^2\Lambda^2Q^2-2a\Lambda Q-\frac14,\\
-\Delta_g^2\Phi-\frac12\Delta_gF
={}&
-a\Lambda\Delta_g^2\varphi-2a\Lambda\Delta_gQ.
\end{align*}
The remaining contribution is $-\frac12\tr_g(\partial_tg)V$.

Let $\mathcal J_j$ denote the terms of degree $j$ in $a$, and let $\mathcal J=\mathcal J_3+\mathcal J_2+\mathcal J_1+\mathcal J_0$. Grouping the preceding terms according to their degree in $a$, we obtain
\begin{equation}\label{eq:carleman-J-lower-bound}
\begin{aligned}
&\int e^{2\Phi}|\partial_tu+\Delta_gu|^2\,\dd\mu_g\dd t\\
&\quad\ge
\int\Bigg\{
\Bigg[
4a\Lambda(\nabla_g^2\varphi+Qg)
+
\left(1+\frac12\tr_g(\partial_tg)\right)g
-\partial_tg
\Bigg](\nabla_gY,\nabla_gY)\\
&\hspace{42mm}
+
(\mathcal J_3+\mathcal J_2+\mathcal J_1+\mathcal J_0)Y^2
\Bigg\}\,\dd\mu_g\dd t.
\end{aligned}
\end{equation}

The four zero-order terms are
\begin{equation}\label{eq:J3}
\mathcal J_3
=
4a^3\Lambda^3
\left[
\nabla_g^2\varphi(\nabla_g\varphi,\nabla_g\varphi)
-Q|\nabla_g\varphi|^2
\right],
\end{equation}
\begin{equation}\label{eq:J2}
\begin{aligned}
\mathcal J_2
={}&
4a^2\Lambda(-\Lambda')
\bigl(|\nabla_g\varphi|^2-\varphi Q\bigr)
-a^2\Lambda^2\bigl(|\nabla_g\varphi|^2+4Q^2\bigr)\\
&+
a^2\Lambda^2
\left[
\partial_tg-\frac12\tr_g(\partial_tg)g
\right](\nabla_g\varphi,\nabla_g\varphi),
\end{aligned}
\end{equation}
\begin{equation}\label{eq:J1}
\begin{aligned}
\mathcal J_1
={}&
a(\Lambda''+\Lambda')\varphi
-a\Lambda\Delta_g^2\varphi
-2a\Lambda\Delta_gQ
+a(8t+2)\Lambda Q\\
&+
\frac a2\tr_g(\partial_tg)\Lambda'\varphi,
\end{aligned}
\end{equation}
and
\begin{equation}\label{eq:J0}
\mathcal J_0
=
\frac{11}{4}+2t+s\tr_g(\partial_tg).
\end{equation}

We estimate these terms in descending order. By \eqref{eq:transfer_estimates}, $\mathcal J_3\ge ca^3\Lambda^3\rho^{1/2}\ge0$.

For the quadratic term, retain one copy of $|\nabla_g\Phi|^2=a^2\Lambda^2|\nabla_g\varphi|^2$. Using \eqref{eq:time_bounds}, \eqref{eq:transfer_estimates}, and $|\partial_tg|\le C\rho^{-2}$, we obtain
\[
\begin{aligned}
\mathcal J_2-|\nabla_g\Phi|^2
\ge{}&
4a^2\Lambda(-\Lambda')
\left[
\left(1-\frac{\Lambda}{2(-\Lambda')}\right)
|\nabla_g\varphi|^2-\varphi Q
\right]\\
&-Ca^2\Lambda^2\rho^{-1}
\ge
ca^2\Lambda(-\Lambda')\rho.
\end{aligned}
\]
Here the expression in brackets is bounded below by $c\rho$, and the error is absorbed after increasing $R$, since $\Lambda\le C(-\Lambda')$.

For the linear term, $a(8t+2)\Lambda Q\ge0$. The remaining estimates in \eqref{eq:transfer_estimates}, together with $|\tr_g(\partial_tg)|\le C\rho^{-2}$, give
\[
\mathcal J_1
\ge
a(\Lambda''+\Lambda')\varphi
-Ca\Lambda
-Ca\rho^{-2}(-\Lambda')\varphi.
\]
By \eqref{eq:time_bounds}, the last term is absorbed by the first after increasing $R$. Hence $\mathcal J_1\ge ca\varphi-Ca\Lambda$.

Finally, $\mathcal J_0\ge c-C\rho^{-2}\ge c$ for sufficiently large $R$. Therefore
\[
\mathcal J
\ge
|\nabla_g\Phi|^2
+ca^2\Lambda(-\Lambda')\rho
+ca\varphi-Ca\Lambda+c.
\]
Choose $a_0$ sufficiently large so that the quadratic term absorbs $Ca\Lambda$ and leaves a positive multiple of $a\Lambda$. It follows that
\[
\mathcal J
\ge
|\nabla_g\Phi|^2
+c\bigl(1+a(\Lambda+\varphi)\bigr).
\]

Together with $\mathcal G\ge cg$, this gives
\[
\int e^{2\Phi}|\partial_tu+\Delta_gu|^2
\ge
c\int
\left[
|\nabla_gY|^2
+|\nabla_g\Phi|^2Y^2
+\bigl(1+a(\Lambda+\varphi)\bigr)Y^2
\right].
\]
Since $e^\Phi\nabla_gu=\nabla_gY-Y\nabla_g\Phi$, we have
\[
e^{2\Phi}|\nabla_gu|^2
\le
C\left(
|\nabla_gY|^2+|\nabla_g\Phi|^2Y^2
\right).
\]
This proves \eqref{eq:carleman}.
\end{proof}

\section{Proof of Theorem~\ref{thm:main}} \label{sec:proof}
In this section, let $\Sigma^{\circ}$ denote the interior of $\Sigma.$

We now apply the Carleman estimate to a cutoff of $w$ supported away from the angular side, infinity, and $t=0$. The Gaussian decay controls the radial and initial-time errors, while the  cutoff is handled by the Carleman weight.

\begin{proof}[Proof of Theorem~\ref{thm:main}]
Choose $R$ sufficiently large so that Proposition \ref{prop:carleman} holds and the lower-order terms below can be absorbed. Fix $B$ so that $\Lambda\varphi<B/4$ on $\{\rho\le3R\}$. Let $\chi_R=0$ on $\{\rho\le2R\}$ and $\chi_R=1$ on $\{\rho\ge3R\}$, and choose $\chi\in C^\infty(\R)$ non-decreasing with $\chi=0 $ on $(-\infty,1/4],$ $\chi=1 $ on $ [1/2,\infty).$ 
Let $\eta=\chi_R(\rho)\chi\left(\frac{\Lambda\varphi}{B}\right).$ Since the second factor vanishes on $\operatorname{supp}\nabla\chi_R$, all derivatives of $\eta$ are supported in $\mathcal B=\{B/4<\Lambda\varphi<B/2\}.$ On $\mathcal B$, $e^{2\Phi-2aB}\le Ce^{-aB}.$

For $L>4R$, let $\zeta_L(\rho)=\zeta(\rho/L)$, where $\zeta=1$ on $(-\infty,1]$ and $\zeta=0$ on $[2,\infty)$. Then $|\nabla_g\zeta_L|\le C/L$ and $|\Delta_g\zeta_L|\le C/L^2$. Let $\theta_\delta=0$ for $t\le\delta$, $\theta_\delta=1$ for $t\ge2\delta$, and $|\theta_\delta'|\le C/\delta$. Let
\[
u_{L,\delta}=\theta_\delta\zeta_L\eta w.
\]
On its support, $\rho\le2L$ and $\Lambda\varphi>B/4$. Since
\[
\varphi=\rho^{3/2}\left[\left(\frac r\rho\right)^{3/2}-\frac18\right],
\]
these conditions give a positive separation from $r=\rho/4$ for fixed $L$. Moreover, $\Lambda(t)\varphi\to0$ as $t\uparrow1/2$, uniformly on $\{\rho\le2L\}$. Hence the cutoff vanishes near $t=1/2$, and $u_{L,\delta}\in C_c^\infty(\mathcal D_R\times(0,1/2)).$

On $\mathcal B$, 
\[
\nabla_g\eta =\chi_R\chi'\left(\frac{\Lambda\varphi}{B}\right) \frac{\Lambda}{B}\nabla_g\varphi,
\]
and
\[
\begin{aligned}
(\partial_t+\Delta_g)\eta
={}&\frac{\chi_R}{B}\chi'
\left(\frac{\Lambda\varphi}{B}\right)
(\Lambda'\varphi+\Lambda\Delta_g\varphi)
+\frac{\chi_R\Lambda^2}{B^2}\chi''
\left(\frac{\Lambda\varphi}{B}\right)|\nabla_g\varphi|^2.
\end{aligned}
\]
Hence, by Lemma~\ref{lem:transfer},
\[
|\nabla_g\eta|\le C_B\rho^{1/2},\qquad
|(\partial_t+\Delta_g)\eta|\le C_B(1+\rho^{3/2}).
\]

The product rule gives
\[
\begin{aligned}
(\partial_t+\Delta_g)u_{L,\delta}
={}&\theta_\delta\zeta_L\eta(\partial_t+\Delta_g)w
+2\theta_\delta\zeta_L\langle\nabla_g\eta,\nabla_gw\rangle\\
&+\theta_\delta\zeta_Lw(\partial_t+\Delta_g)\eta
+2\theta_\delta\eta\langle\nabla_g\zeta_L,\nabla_gw\rangle\\
&+\theta_\delta\eta w\Delta_g\zeta_L
+2\theta_\delta w\langle\nabla_g\zeta_L,\nabla_g\eta\rangle
+\theta_\delta'\zeta_L\eta w.
\end{aligned}
\]
Since
\[
\theta_\delta\zeta_L\eta\,\nabla_gw
=\nabla_gu_{L,\delta}
-w\nabla_g(\theta_\delta\zeta_L\eta),
\]
Proposition \ref{prop:difference} gives
\[
|\theta_\delta\zeta_L\eta(\partial_t+\Delta_g)w|^2
\le CR^{-4}u_{L,\delta}^2
+CR^{-6}|\nabla_gu_{L,\delta}|^2
+C\rho^{-6}w^2|\nabla_g(\theta_\delta\zeta_L\eta)|^2.
\]
The first two terms are absorbed after increasing $R$.

Apply Proposition \ref{prop:carleman} to $u_{L,\delta}$ with weight $e^{2\Phi-2aB}$. The remaining terms are supported in $\mathcal B$, $\{\delta<t<2\delta\}$, or $\{L<\rho<2L\}$.

On $\mathcal B$, Proposition~\ref{lem:gaussian} and the sectorial volume bound give
\[
C_Be^{-aB}
\int_0^{1/2}\int_{\mathcal D_R}
(1+\rho^3)(w^2+|\nabla_gw|^2)\,d\mu_gdt
\le C_Be^{-aB}.
\]
For fixed $a$ and $L$,
\[
E_\delta
\le C_{a,L}\delta^{-2}
\int_\delta^{2\delta}e^{-cR^2/t}\,dt
\le C_{a,L}\delta^{-1}e^{-cR^2/\delta}\to0,
\]
while, for fixed $a$,
\[
E_L\le C_aL^N e^{CaL^{3/2}-cL^2}\to0.
\]
Letting first $\delta\downarrow0$ and then $L\to\infty$, Fatou's lemma gives
\begin{equation}\label{eq:limit_inequality}
\begin{aligned}
&\int_0^{1/2}\int_{\mathcal D_R}e^{2\Phi-2aB}
\left(
|\nabla_g(\eta w)|^2
+\bigl(1+a(\Lambda+\varphi)\bigr)(\eta w)^2
\right)\,d\mu_gdt \le C_Be^{-aB}.
\end{aligned}
\end{equation}

Choose $0<t_1<t_2<1/2$. Since $\Lambda(t_2)>0$, take $\rho_*>3R$ so large that $\Lambda(t_2)\varphi\ge2B$ whenever $\rho_*\le\rho\le\rho_*+1$ and $r\ge\rho/2$. Let
\[
E=\{\rho_*\le\rho\le\rho_*+1,\ r\ge\rho/2,\ t_1\le t\le t_2\}.
\]
Then $\eta=1$ on $E$ and $\Phi-aB\ge aB$. Hence
\[
\int_Ew^2\,d\mu_gdt\le C_Be^{-3aB}.
\]
Letting $a\to\infty$ gives $w=0$ on the interior of $E$. Thus, for any $t_0\in(t_1,t_2)$, the two rescaled ends agree on a nonempty open set. Scaling back by $t_0^{-1/2}$, we conclude that $\Sigma$ and $\widetilde\Sigma$ agree on a nonempty open set.

Let $\mathcal Z\subset\Sigma^\circ$ 
be the set of points near which the two ends coincide. Then $\mathcal Z$ is nonempty and open. To show that it is closed, let $P_j\in\mathcal Z$ with $P_j\to P\in\Sigma^\circ$. By properness and smoothness, $P\in\widetilde\Sigma^\circ$ and $T_P\Sigma=T_P\widetilde\Sigma$. Near $P$, write the two ends as graphs $F_i(\xi)=P+\xi+v_i(\xi)\nu_*$, $i=1,2$, over their common tangent plane. Let
\[
a^{j\ell}(p)=\delta^{j\ell}-\frac{p_jp_\ell}{1+|p|^2}.
\]
The self-shrinker equation becomes
\[
a^{j\ell}(Dv_i)\partial_{j\ell}v_i
-\frac12(P^T+\xi)\cdot Dv_i
+\frac12v_i+\frac12P\cdot\nu_*=0,
\]
where $P^T$ is the tangential projection of $P$. Hence $v=v_2-v_1$ satisfies
\begin{equation}\label{eq:elliptic_diff}
a^{j\ell}(Dv_2)\partial_{j\ell}v
+b^q\partial_qv+\frac12v=0,
\end{equation}
where $b^q$ is smooth. Since
\[
\frac{|\zeta|^2}{1+|Dv_2|^2}
\le a^{j\ell}(Dv_2)\zeta_j\zeta_\ell
\le|\zeta|^2,
\]
the equation is locally uniformly elliptic. As $v$ vanishes near $P_j$ for sufficiently large $j$, elliptic unique continuation \cite{GL} gives $v=0$ near $P$. Thus $\mathcal Z$ is closed in $\Sigma^\circ$. By connectedness, $\mathcal Z=\Sigma^\circ$, and interchanging the two ends yields $\Sigma=\widetilde\Sigma$.
\end{proof}

\section{Extension to higher codimension}\label{extension}
In this section we prove Theorem~\ref{thm:2} and Corollary~\ref{thm3}. The scalar graph functions are now replaced by sections of the normal bundle of $\K$. The detailed normal-graph calculations and vector estimates are given in Appendix~\ref{app:higher-codimension}.

Let $\C^m \subset \R^{m+n}$ be as in Theorem~\ref{thm:2}. We use the same functions $r, \rho, \varphi, Q, \Lambda, \Phi$ as above. The induced metric on $\K$ is still
\begin{align*}
g_0&=dr^2+r^2g_\Gamma+\sum_{\mu=1}^k (dz_\mu)^2.
\end{align*}
Let $N\K$ denote the normal bundle of $\K$, let $A_0$ denote its  normal-vector-valued second fundamental form, and let $\nabla^{\perp}$ denote its normal connection. Covariant derivatives of normal-bundle-valued tensor use $\nabla^{\perp}$ on their normal indices and $\nabla_0$ on their tangent indices. The rescaled ends have the normal graph representations
\[
F(X,t)=X+U(X,t),\qquad \widetilde F(X,t)=X+\widetilde U(X,t),\qquad 0<t\le1,
\]
where $U(X,t),\widetilde U(X,t)\in N_X\K$. For each fixed $t$, the maps $F(\cdot,t)$ and $\widetilde F(\cdot,t)$ describe the embedded normal graphs over the sector specified below. We write $w=\widetilde U-U$ and $g(t)=g_{U(t)}$. When $w$ is viewed as a function into $\R^{m+k+n}$, the operators $\nabla_0$, $\nabla_g$, and $\Delta_g$ act on its fixed Euclidean components. The operator $\Delta_g^\perp$ uses the fixed connection $\nabla^\perp$ and the metric $g(t)$ on the domain.

We first establish the normal-bundle analogues of the graphical and metric estimates in Section~\ref{sec:geometry}, together with a differential inequality for the difference $w$.

\begin{proposition}\label{prop normal}
There is $R>0$ such that the normal graphs are defined over $\{\rho>R,\ r>\rho/16\}$ for $0<t\le1$. On $\{\rho>R,\ r>\rho/8\}\times(0,1]$, for every fixed nonnegative integer $j$, we have
\begin{align}
|(\nabla^\perp)^jU|+|(\nabla^\perp)^j\widetilde U|&\le C_jt\rho^{-1-j},\label{eq:normal-graph-bounds}\\
|(\nabla^\perp)^j\partial_tU|+|(\nabla^\perp)^j\partial_t\widetilde U|&\le C_j\rho^{-1-j}.\label{eq:normal-time-bounds}
\end{align}
Moreover, the induced metric satisfies
\begin{equation}\label{eq:normal-metric}
\begin{aligned}
\sum_{j=0}^3\rho^j|\nabla_0^j(g-g_0)|&\le Ct\rho^{-2},\\
|\partial_tg|+\rho|\nabla_0\partial_tg|&\le C\rho^{-2},\qquad \tfrac12g_0\le g\le2g_0.
\end{aligned}
\end{equation}
The difference $w$ satisfies
\begin{equation}\label{eq:normal-rough-difference}
|\partial_tw+\Delta_g^\perp w|\le C\left(\rho^{-2}|w|+t\rho^{-3}|\nabla^\perp w|_g\right).
\end{equation}
The metric norms in \eqref{eq:normal-metric} are taken with respect to $g_0$. The constants depend on $m,k,n,\C,M_1,M_2$, and $C_j$ may also depend on $j$.
\end{proposition}
The proof is given in Appendices~\ref{app:normal-calculations}--\ref{app:normal-difference}.

To apply the scalar estimates from Sections~\ref{sec:gaussian} and~\ref{sec:carleman}, we next view $w$ in fixed ambient Euclidean components. The following proposition compares the normal rough Laplacian with the componentwise Laplacian and gives the corresponding differential inequality.

\begin{proposition}\label{prop:ambient-difference}
Let $\Pi(X)$ denote the orthogonal projection onto $N_X\K$. On $\{\rho>R,\ r>\rho/8\}\times(0,1]$, we have
\begin{equation}\label{eq:normal-projection}
|\nabla_g\Pi|\le C\rho^{-1},
\qquad
|\Delta_g\Pi|\le C\rho^{-2}.
\end{equation}
Moreover,
\begin{equation}\label{eq:normal-laplacian-comparison}
|\Delta_gw-\Delta_g^\perp w|
\le C\left(\rho^{-2}|w|+\rho^{-1}|\nabla_gw|\right),
\end{equation}
and consequently,
\begin{equation}\label{eq:diff_pde-normal}
|\partial_tw+\Delta_gw|
\le C\left(\rho^{-2}|w|+\rho^{-1}|\nabla_gw|\right).
\end{equation}
For every fixed nonnegative integer $j$, we also have
\begin{equation}\label{eq:normal-ambient-bounds}
  |\nabla_0^jw|\le C_jt\rho^{-1-j},
\qquad
|\nabla_0^j\partial_tw|\le C_j\rho^{-1-j},  
\end{equation}
and
\begin{equation}\label{eq:w_scaling-normal}
w(X,t)=\sqrt t\,w(X/\sqrt t,1).
\end{equation}
Here $C$ depends on $m,k,n,\C,M_1,M_2$, and $C_j$ also depend on $j$.
\end{proposition}

The proof is given in Appendix~\ref{app:normal-components}.

Using the ambient differential inequality \eqref{eq:diff_pde-normal} and the scaling identity \eqref{eq:w_scaling-normal}, we obtain Gaussian decay by the same localization and scaling argument as in Section~\ref{sec:gaussian}.

\begin{proposition}\label{lem:gaussian-normal}
Under the hypotheses of Proposition~\ref{prop normal}, there are $R,C,c>0$ such that
\begin{equation}\label{eq:normal-gaussian}
|w(X,t)|+|\nabla_{g(t)}w(X,t)|\le C\exp(-c\rho^2/t)
\end{equation}
for $\rho>R$, $r\ge\rho/4$, and $0<t\le1$. The constants depend on $m,k,n,\C,M_1,M_2$.
\end{proposition}
The proof  is given in Appendix~\ref{app:gaussian}.

The metric bounds \eqref{eq:normal-metric} also allow us to use the Carleman estimate from Section~\ref{sec:carleman} with the same weight function. Applying the scalar estimate to each fixed Euclidean component and summing yields the following inequality.

\begin{proposition}\label{prop:carleman-normal}
Suppose that $g(t)$ satisfies \eqref{eq:normal-metric}. There are $R,a_0,C>0$ such that, for every $a\ge a_0$ and every $u\in C_c^\infty(\mathcal D_R\times(0,1/2);\R^{m+k+n})$,
\begin{equation}\label{eq:normal-carleman}
\begin{aligned}
&\int_0^{1/2}\!\int_{\mathcal D_R}e^{2\Phi}\left[|\nabla_gu|^2+\bigl(1+a(\Lambda+\varphi)\bigr)|u|^2\right]\dd\mu_g\dd t\\
&\qquad\le C\int_0^{1/2}\!\int_{\mathcal D_R}e^{2\Phi}|\partial_tu+\Delta_gu|^2\dd\mu_g\dd t.
\end{aligned}
\end{equation}
The constants depend on $m,k,n,\C$ and the constants in \eqref{eq:normal-metric}, and are independent of $a$ and $u$.
\end{proposition}
The proof is given in Appendix~\ref{app:carleman}.

\begin{proof}[Proof of Theorem~\ref{thm:2}]
Propositions~\ref{prop normal}--\ref{prop:carleman-normal}, together with the cutoff argument in Appendix~\ref{app:carleman}, show that $w$ vanishes on a nonempty open set at a positive time. Scaling back, $M_1$ and $M_2$ agree on a nonempty open set. The same open-and-closed argument as in the proof of Theorem~\ref{thm:main}, together with Lemma~\ref{lem:local-shrinker-continuation}, yields $M_1=M_2$.
\end{proof}

\begin{proof}[Proof of Corollary~\ref{thm3}]
    The product manifold $M_2 \times \R^k$ has induced metric $g_{M_2}+\sum^k_{\mu=1}d\,z^2_{\mu}.$ Its mean curvature vector and  the normal projection of its position vector satisfy
    \[
    \mathbf{H}_{M_2\times\R^k} 
    = (\mathbf{H}_{M_2}, 0)
    =- \frac{1}{2}(X_{M_2}^{\perp}, 0)
    =- \frac{1}{2}(X_{M_2}, z)^{\perp}.
    \]
Thus $M_2\times\R^k$ is a connected, complete, properly embedded self-shrinker without boundary. Moreover, 
\[
\lambda(M_2\times \R^k)
= (\lambda M_2) \times \R^k
\longrightarrow \C \times \R^k
\]
locally smoothly away from the axis. The same local comparison and continuation argument as in the proof of Theorem~\ref{thm:2} gives $M_1=M_2\times\R^k$.
    
\end{proof}

\appendix
\section{ Normal-bundle calculation} \label{app:higher-codimension}

In this appendix, we derive the normal graph equations and the difference estimates used in Section~\ref{extension}. For the Gaussian decay, Carleman estimate, and cutoff argument, we refer to the codimension-one proofs and explain only the changes needed for normal sections. The final subsection verifies the hypotheses of elliptic unique continuation for the local graph difference.

Throughout the appendix, $\nabla^\perp$ is the fixed normal connection on $N\K$, and $\nabla_0$ is the Levi-Civita connection of $g_0$. Derivatives of tensors with normal indices use both connections as described in Section~\ref{extension}. For ambient-vector-valued functions, derivatives act on fixed Euclidean components. Unless indicated otherwise, norms are taken with respect to $g_0$ and the Euclidean inner product. Constants may depend on $m,k,n,\C,M_1,M_2$, the fixed sector, and the indicated derivative order.

\subsection{The normal graph equation}\label{app:normal-calculations}
The proofs of Lemmas~\ref{lem:sector} and~\ref{lem:graphs} apply with normal covariant derivatives. They give the same sector estimates, graph domain, coarse derivative bounds, uniform smallness at infinity, and local convergence of all spatial derivatives to zero as $t\downarrow0$. We record only the changes in the graphical equation~\eqref{eq:graph_equation}. For the normal graph $F(X,t)=X+U(X,t)$, define the shape operator $A_U^\sharp$ and the tangent endomorphism $S_U$ by
\[
g_0(A_U^\sharp V,W)=\langle A_0(V,W),U\rangle,\qquad S_U=I-A_U^\sharp.
\]
In a local $g_0$-orthonormal frame, the Weingarten formula gives $dF(e_i)=S_Ue_i+\nabla_i^\perp U$, and hence
\begin{equation}\label{eq:app-normal-metric-formula}
(g_U)_{ij}=g_0(S_Ue_i,S_Ue_j)+\langle\nabla_i^\perp U,\nabla_j^\perp U\rangle.
\end{equation}
Let $T_U$ be the tangent-vector-valued tensor
\[
(T_U)_{ij}=(\nabla_{0,i}A_U^\sharp)e_j+A_{\nabla_j^\perp U}^\sharp e_i.
\]
The differentiation used for \eqref{eq:graph_equation} gives $(\nabla_0^2F)_{ij}^{\top}=-(T_U)_{ij}$, where $\top$ denotes projection onto $T\K$. Since the Christoffel correction is tangent to the graph, backward mean curvature flow implies that $\partial_tF+g_U^{ij}(\nabla_0^2F)_{ij}$ is graph tangent. Such vectors have the form $S_UV+\nabla_V^\perp U$, and $S_U$ is invertible for large $R$. Comparing tangent and normal components therefore gives
\begin{equation}\label{eq:app-normal-graph-equation}
\begin{aligned}
0=\partial_tU+g_U^{ij}\bigl[&((\nabla^\perp)^2U)_{ij}+A_0(e_i,e_j)-A_0(e_i,A_U^\sharp e_j)
+\nabla_{S_U^{-1}(T_U)_{ij}}^\perp U\bigr].
\end{aligned}
\end{equation}

\subsection{Estimates and the difference equation}\label{app:normal-difference}
The proof of Lemma~\ref{lem:fine} applies with the fixed normal connection and gives \eqref{eq:normal-graph-bounds}--\eqref{eq:normal-time-bounds}. That argument uses only the product rule, differentiation of inverses, and integration from $t=0$; it does not commute spatial covariant derivatives. The metric estimates \eqref{eq:normal-metric} follow from \eqref{eq:app-normal-metric-formula} by the same differentiation as in Proposition~\ref{prop:difference}. In particular,
\begin{equation}\label{eq:app-normal-connection-bounds}
|\nabla_0^j(\Gamma_g-\Gamma_0)|\le Ct\rho^{-3-j},\qquad 0\le j\le2.
\end{equation}
For $w=\widetilde U-U$, subtract the two normal graph equations as in Proposition~\ref{prop:difference}. The matrix inverse identity used there remains valid for the shape endomorphisms:
\[
S_{\widetilde U}^{-1}-S_U^{-1}=S_{\widetilde U}^{-1}A_w^\sharp S_U^{-1}.
\]
Keeping this order of the factors, the same product estimates, now with normal inner products and tensor norms, give
\[
\left|\partial_tw+g^{ij}((\nabla^\perp)^2w)_{ij}\right|\le C\bigl(\rho^{-2}|w|+t\rho^{-3}|\nabla^\perp w|\bigr).
\]
The normal rough Laplacian for the fixed normal connection and the metric $g(t)$ is
\begin{equation}\label{eq:app-normal-laplacian}
\Delta_g^\perp w=g^{ij}((\nabla^\perp)^2w)_{ij}-g^{ij}(\Gamma_g-\Gamma_0)^\ell_{ij}\nabla_\ell^\perp w.
\end{equation}
The connection estimate and metric comparison therefore give \eqref{eq:normal-rough-difference}, completing the proof of Proposition~\ref{prop normal}.

\subsection{Reduction to fixed Euclidean components}\label{app:normal-components} We prove Proposition~\ref{prop:ambient-difference}. The orthogonal projection $\Pi(X)$ onto $N_X\K$ is a globally defined smooth matrix on $\K$. Since $N_{\lambda X}\K=N_{ X}\K$ as subspaces of $\R^{m+k+n}$, it is homogeneous of degree zero. Compactness of normalized closed sectors and scaling give
\[
|\nabla_0\Pi(X)| \leq C \rho^{-1},
\qquad
|\nabla^2_0\Pi(X)| \leq C \rho^{-2}
\]
on $r> \rho/8.$ The metric and connection estimate give 
\[
|\nabla_g \Pi(X)|\leq C\rho^{-1},
\qquad
|\Delta_g \Pi(X)|\leq C\rho^{-2}.
\]

In local coordinates, $w=\Pi w$ and $\nabla^{\perp}_i w= \Pi( \partial_i w)$. Expanding the normal Laplacian gives 
\begin{align*}
  \Delta^{\perp}_g w
&=g^{ij}[\Pi \partial_i(\Pi \partial_jw)-\Gamma^l_{g,ij} \Pi \partial_l w]\\
&= \Pi \Delta_g w+ g^{ij}\Pi(\partial_i \Pi)\partial_j w.
\end{align*}
On the other hand, applying $\Delta_g$ to $w=\Pi w$ gives
\begin{align*}
    \Delta_g w
    &= (\Delta_g\Pi) w+ 2 g^{ij} (\partial_i \Pi)\partial_j w + \Pi \Delta_g w. 
\end{align*}
Combining the last two identities, we obtain
\begin{equation}\label{eq:app-projection-comparison}
    \Delta_g w- \Delta_g ^{\perp} w
    = (\Delta_g \Pi) w +  g^{ij} (2I-\Pi)(\partial_i \Pi) \partial_j w.
\end{equation}
Since $\Pi$ is an orthogonal projection, then $|\Pi(X)|\leq 1$, $|2I- \Pi(X)|\leq 2,$ we get
\[
 |\Delta_g w- \Delta_g ^{\perp} w| 
    \leq C(\rho ^{-2}|w|+ \rho^{-1}|\nabla_g w|).
\]
Also, $|\nabla^\perp w|_g\le|\nabla_gw|$. Combining the preceding inequality with \eqref{eq:normal-rough-difference} proves \eqref{eq:diff_pde-normal}.
The ambient derivative bounds follow from the Weingarten formula
\[
\overline{\nabla}_V w
= - A^\#_w V+ \nabla^{\perp}_V w.
\]
Further differentiation produces derivatives of $A_0$ and normal covariant derivatives of $w$. Every derivative introduces one additional factor $\rho^{-1}$. Induction using the sector estimates from Lemma~\ref{lem:sector}, \eqref{eq:normal-graph-bounds}, and \eqref{eq:normal-time-bounds} proves \eqref{eq:normal-ambient-bounds} for all fixed $j$.

Finally, let $Y = X/\sqrt{t}$. The normal spaces at $X$ and $Y$ coincide as subspaces of $\R^{m+k+n}$, and 
\[
\sqrt{t}[Y+ U(Y,1)]= X+ \sqrt{t}  U(X/\sqrt{t}, 1).
\]
Uniqueness of normal projection gives $U(X, t)=\sqrt{t}U(X/\sqrt t,1)$, and the same identity holds for $\widetilde U$. Subtracting the two identities proves \ref{eq:w_scaling-normal}. Differentiation under rescaling gives
\[
|\nabla^j_0 w(X,t)|
= t^{(1-j)/2}|\nabla_0^j w(X/\sqrt{t}, 1)|.
\]
This completes the proof of Proposition~\ref{prop:ambient-difference}.

\subsection{Gaussian decay}\label{app:gaussian} Lemma~\ref{lem:local_l2} extends to functions with values in $\R^{m+k+n}$. Apply \eqref{eq:local_carleman} to each fixed Euclidean component of $\widehat u=\eta\chi_\alpha u$, and sum before using the vector differential inequality. The sum of the supremum remainders is bounded by a fixed multiple of the supremum of $\int(|\widehat u|^2+|D\widehat u|^2)\dd x$. The absorption, cutoff estimates, and approximation at $t=0$ are then those in the scalar lemma, with constants also depending on the number of components.

We now follow the proof of Proposition~\ref{lem:gaussian}. In its tangent-plane chart centered at $X_0,$ with $\rho_0=|X_0|$, the metric bounds \eqref{eq:normal-metric} and the ambient differential inequality \eqref{eq:diff_pde-normal} give 
\[
|\partial_t w+ \partial_i(g^{ij}\partial_j w)|
\leq C (\rho_0^{-2}|w|+\rho_0^{-1}|Dw|).
\]
This is the local divergence-form bound used in that proof; its parabolic rescaling makes the first-order coefficient uniformly bounded. The initial vanishing and derivative bounds needed for the local estimate and interpolation follow from \eqref{eq:normal-ambient-bounds}. Applying the vector version of Lemma~\ref{lem:local_l2}, followed by the same interpolation and scaling argument using \eqref{eq:w_scaling-normal}, proves \eqref{eq:normal-gaussian}, after adjusting the final radius as in the scalar proof.

\subsection{Carleman estimate}\label{app:carleman} The proof of Proposition~\ref{prop:carleman} only uses the product metric and the metric estimates established above, and therefore applies in arbitrary codimension. Applying the scalar estimate to each fixed Euclidean component and summing gives Proposition~\ref{prop:carleman-normal}. We use the same cutoffs as in Section~\ref{sec:proof} and write $f_{L,\delta}=\theta_\delta \zeta_L \eta$ and $u_{L,\delta}= f_{L,\delta}w.$ By \eqref{eq:diff_pde-normal},
\[
|f_{L,\delta}(\partial_t+\Delta_g) w|^2
\leq CR^{-4}|u_{L,\delta}|^2 
+ CR^{-2}|\nabla_g u_{L,\delta}|^2+C\rho^{-2}|w|^2|\nabla_gf_{L,\delta}|^2.
\]
The first two terms are absorbed into \eqref{eq:normal-carleman} by choosing $R $ sufficiently large. Moreover, for $\mathcal B=\{B/4\le\Lambda\varphi\le B/2\}$, $\mathcal I=\{L\leq \rho\leq 2L\}$, we get
\[
\rho^{-2}|\nabla_gf_{L,\delta}|^2
\leq C_B\rho^{-1} \mathbf 1_{\mathcal B} + CL^{-4} \mathbf 1_\mathcal I.
\]
These errors are controlled by the Gaussian decay estimate\eqref{eq:normal-gaussian}. Using \eqref{eq:normal-carleman}, we follow the limiting argument in Section~\ref{sec:proof}: first let $\delta\downarrow0$, then $L\to\infty$ for fixed $a$, and finally $a\to\infty$. It follows that $w$ vanishes on a nonempty open set at a positive time.

\subsection{Unique continuation for elliptic systems}\label{app:UCT}
\begin{lemma}\label{lem:local-shrinker-continuation}
Let $f_1,f_2:\Omega\to\R^n$ be smooth functions whose graphs over the same affine plane are self-shrinkers, where $\Omega\subset\R^{m+k}$ is a connected open set. If $f_1=f_2$ on a nonempty open subset of $\Omega$, then $f_1=f_2$ throughout $\Omega$.
\end{lemma}
\begin{proof}
As in the hypersurface case, the difference $v=f_2-f_1$ satisfies a linear elliptic system of the form
\[
a^{ij}(Df_2)\partial_{ij}v^\alpha+B^{\ell\alpha}{}_\beta\partial_\ell v^\beta+\frac12v^\alpha=0,
\]
where $a^{ij}(Df_2)$ is the inverse of the induced metric and $B^{\ell\alpha}{}_\beta$ are smooth coefficients. The principal part is the same locally uniformly elliptic scalar operator on every component, with only lower-order coupling. The conclusion follows from Aronszajn's unique continuation theorem for elliptic systems \cite[Remark~3, p.~248]{Aronszajn}.
\end{proof}

\end{document}